\documentclass[11pt,english]{extarticle}
\usepackage{lmodern}
\usepackage[T1]{fontenc}
\usepackage[letterpaper]{geometry}
\usepackage{color}
\usepackage{babel}
\usepackage{array}
\usepackage{float}
\usepackage{booktabs}
\usepackage{multirow}
\usepackage{amsmath}
\usepackage{amsthm}
\usepackage{amssymb}
\usepackage{graphicx}
\usepackage{setspace}
\usepackage[authoryear]{natbib}
\usepackage[unicode=true,
 bookmarks=true,bookmarksnumbered=true,bookmarksopen=false,
 breaklinks=false,pdfborder={0 0 0},pdfborderstyle={},backref=false,colorlinks=true]
 {hyperref}

\makeatletter

\newcommand{\noun}[1]{\textsc{#1}}
\providecommand{\tabularnewline}{\\}
\floatstyle{ruled}
\newfloat{algorithm}{tbp}{loa}
\providecommand{\algorithmname}{Algorithm}
\floatname{algorithm}{\protect\algorithmname}

\theoremstyle{definition}
\newtheorem{example}{\protect\examplename}
\theoremstyle{remark}
\newtheorem{rem}{\protect\remarkname}
\theoremstyle{definition}
\newtheorem{defn}{\protect\definitionname}
\theoremstyle{plain}
\newtheorem{prop}{\protect\propositionname}
\theoremstyle{plain}
\newtheorem{lem}{\protect\lemmaname}

\@ifundefined{date}{}{\date{}}
\usepackage{xcolor}

\hypersetup{
	colorlinks,
	linkcolor={red!60!black},
	citecolor={blue!60!black},
	urlcolor={blue!80!black}
}

\usepackage{algorithm,algpseudocode}

\makeatother

\providecommand{\definitionname}{Definition}
\providecommand{\examplename}{Example}
\providecommand{\lemmaname}{Lemma}
\providecommand{\propositionname}{Proposition}
\providecommand{\remarkname}{Remark}

\begin{document}
\title{A Flexible Multi-Facility Capacity Expansion Problem \\
	 with Risk Aversion}

\maketitle
\vspace{-2cm}
\singlespacing
\begin{center}
{\large{}Sixiang Zhao}{\large\par}
\par\end{center}
\vspace{-0.7cm}
\begin{center}
{\scriptsize{}Department of Industrial Systems Engineering and Management,
National University of Singapore,  Singapore 117576, \href{mailto:zhaosixiang@u.nus.edu}{zhaosixiang@u.nus.edu}}{\footnotesize\par}
\par\end{center}

\begin{center}
{\large{}William B. Haskell}{\large\par}
\par\end{center}
\vspace{-0.7cm}
\begin{center}
{\scriptsize{}Department of Industrial Systems Engineering and Management,
National University of Singapore, Singapore 117576, \href{mailto:isehwb@nus.edu.sg}{isehwb@nus.edu.sg}}{\footnotesize\par}
\par\end{center}

\begin{center}
{\large{}Michel-Alexandre Cardin}{\large\par}
\par\end{center}
\vspace{-0.7cm}
\begin{center}
{\scriptsize{}Dyson School of Design Engineering, Imperial College
London, London SW7 2AZ, United Kingdom\\ \href{mailto:m.cardin@imperial.ac.uk}{m.cardin@imperial.ac.uk}}{\footnotesize\par}
\par\end{center}

\doublespacing
\begin{abstract}
\emph{This paper studies flexible multi-facility capacity expansion
with risk aversion. In this setting, the decision maker can periodically
expand the capacity of facilities given observations of uncertain
demand. We model this situation as a multi-stage stochastic programming
problem. We express risk aversion in this problem through conditional
value-at-risk (CVaR), and we formulate a mean-CVaR objective. To solve
the multi-stage problem, we optimize over decision rules. In particular,
we approximate the full policy space of the problem with a tractable
family of if-then policies. Subsequently, a decomposition algorithm
is proposed to optimize the decision rule. This algorithm decomposes
the model over scenarios and it updates solutions via the subgradients
of the recourse function. We demonstrate that this algorithm can quickly
converge to high-performance policies. To illustrate the practical
effectiveness of this method, a case study on the waste-to-energy
system in Singapore is presented. These simulation results show that
by adjusting the weight factor of the objective function, decision
makers are able to trade off between a risk-averse policy that has
a higher expected cost but a lower value-at-risk, and a risk-neutral
policy that has a lower expected cost but a higher value-at-risk risk.}

\end{abstract}

\section{Introduction}

The capacity expansion problem aims to determine the capacity plan---that
is, the optimal amount and timing of capacity acquisition---to address
growth of demand. This problem has been widely studied for a variety
of systems, such semiconductor manufacturing \citep{geng_stochastic_2009},
airport facilities \citep{sun_stochastic_2015}, and waste-to-energy
systems \citep{cardin_analyzing_2015}. This problem is challenging
because future demand is uncertain, but capital expenditure for capacity
is high and these investments are usually irreversible (for example,
the capacity of airport facilities, such as highway links and ports,
is hard to decrease once established \citep{sun_stochastic_2015}).
The traditional 'inflexible' method determines the capacity plan at
the beginning before demand is observed, and then this plan is implemented
regardless of the realizations of future demand. However, this method
may suffer unexpected costs if demands fail to grow as expected.

To cope with this issue, increasing attention has been paid to designing
'flexibile' systems in capacity expansion problems (a.k.a. real options
analysis). A flexible system has the ability to adjust its capacity
dynamically as demands are observed. We can then expand the capacity
(i.e. exercise the option) if demand surges, and do nothing when demand
remain steady. In the multi-facility capacity expansion problem (MCEP),
a flexible system has the option not only to adjust capacity, but
also to switch service between facilities. For example, if one facility
runs out of capacity, we can either expand the capacity of this facility,
or allocate excess demand to an adjacent facility. It has been verified
by many industrial case studies that flexibility can improve the expected
lifecycle performance of engineering systems by 10\% to 30\% compared
with inflexible methods \citep{neufville_flexibility_2011,cardin_approach_2017}.

Evaluating the economic performance of a flexible MCEP is a dynamic
optimization problem with uncertainty. This problem needs to find
the optimal expansion policy, which is a mapping from the historical
data to the capacity decisions. However, the optimal policy is usually
hard to characterize for practically-sized MCEP with existing methods.
In particular, dynamic programming (DP)-based algorithms suffer the
curse of dimensionality (in the action space) when the number of facilities
is large. Stochastic programming algorithms can also be inefficient
because of the non-convexity caused by discrete capacity decisions.

To solve MCEPs with discrete capacity, if--then decision rules have
been proposed to approximate the policy space \citep{cardin_approach_2017,cardin_strategic_2017,zhang_flexibility_2017}.
An if-then rule states that if the capacity gap of a facility exceeds
a threshold then its capacity is expanded to a certain level, and
the capacity is unchanged otherwise. In this framework, we want to
find the best if--then decision rule. This type of rule mimics the
behavior of human decision makers so it is intuitive from a managerial
standpoint. In numerical terms, the decision rule-based method not
only provides high-performance solutions for MECPs, it is also more
scalable than DP or stochastic programming \citep{zhao_decision_2018}.

Classical MCEP models suppose that the decision makers are risk-neutral;
in other words, they maximize the expected reward. However, decision
makers often have their own attitudes about risk. We employ the popular
conditional value-at-risk (CVaR) to measure risk in our model. The
motivation for this choice is two-fold. First, CVaR is a coherent
risk measure and has strong decision-theoretic support. Second, CVaR
is convex and enjoys significant computational advantages compared
to other risk-aware objectives \citep{rockafellar_optimization_2000}.

Our present paper investigates a risk-averse MCEP with a mean-CVaR
objective. Our specific contributions are summarized below:
\begin{enumerate}
\item We formulate the objective of the MCEP as a weighted sum of the expected
cost and the CVaR of cost. We solve this risk-averse problem by using
an if--then decision rule to approximate the policy space. To the
best knowledge of the authors, no studies have yet used decision rules
for risk-averse MCEPs.
\item We prove that the value of the risk-averse MCEP is always smaller
than or equal to the risk-neutral MCEP, implying that decision makers
may prefer to pay less for flexibility if they are risk-averse.
\item We optimize over decision rules (with respect to the mean-CVaR objective)
by using a customized decomposition algorithm. This algorithm can
be viewed as an improved version of the one from \citep{zhao_decision_2018}.
We design subgradient cuts to update the parameters of the decision
rule, which is not only more time-efficient but can also provide problem
insight compared with the algorithm in \citep{zhao_decision_2018}.
\end{enumerate}
The remainder of this paper is organized as follows. Section \ref{sec:Literature_Review}
summarizes the relevant literature. In Section \ref{sec:Model-Formulation},
the risk-neutral MCEP model is presented, and then its risk-averse
counterpart is discussed. If--then decision rules are introduced
in Section \ref{sec:Approximation}. In Section \ref{sec:BAC-D_Algorithm},
we present the risk-averse MCEP with decision rules, and solve it
with our subgradient-based decomposition algorithm. A numerical study
on a multi-facility waste-to-energy (WTE) system is elaborated in
Section \ref{sec:Numerical-Study}. Finally, the strengths and limitations
of the proposed method and the opportunities for future research are
summarized. All proofs may be found in the \hyperref[sec:Appendix-A]{Appendix}.

\section{\label{sec:Literature_Review}Literature Review}

\subsection{Flexible Capacity Expansion Problems}

Capacity expansion problems have been widely studied since the seminal
paper by \citep{manne_capacity_1961}. \citet{manne_capacity_1961}
investigated the trade-off between the discount factor and the economies
of scale in capacity expansion problems with deterministic/stochastic
demand. Many variations on this original model have been studied.
Comprehensive reviews may be found in \citep{luss_operations_1982,van_mieghem_commissioned_2003,martinez-costa_review_2014}.

In the framework of real options analysis, capacity expansion decisions
are viewed as a series of options that can be exercised over time
\citep{dixit_investment_1994}. The main advantage of this framework
comes from its 'wait-and-see' nature; capacity decisions can be exercised
or deferred based on the realizations of uncertainty. \citet{eberly_multi-factor_1997}
studied a multi-factor capacity investment problem and characterized
the structure of the optimal policy as an ISD (invest, stay put, and
disinvest) policy. \citet{kouvelis_flexible_2014} studied a flexible
capacity investment problem and investigated the value of the postponed
capacity commitment option for uncertain demand. Our work differs
from these models because we deal with discrete capacity expansion
decisions. In capacity expansion problems with discrete capacity,
\citet{huang_value_2009} derived an analytical bound for the value
of the multi-stage problem compared to the two-stage problem, but
this result requires linear expansion costs. \citet{cardin_analyzing_2015}
and \citet{zhao_decision_2018} studied MCEPs with nonlinear expansion
costs and used power functions to model economies of scale.

\subsection{Risk Measures}

The objectives in the aforementioned papers are all risk-neutral.
To capture the risk preferences of real decision makers, a variety
of risk measures have been reported in the literature including: utility
functions, mean-variance, value-at-risk (VaR), and CVaR. For example,
\citet{hugonnier_real_2007} extended the standard real options analysis
by introducing a utility function to address risk preferences. \citet{birge_option_2000}
incorporated utility functions into a general linear capacity planning
model, and formulated the problem as a multi-stage stochastic programming
problem. Compared to utility functions, CVaR is more intuitive and
easier to specify. Decision makers can express their risk preferences
by directly adjusting the percentile terms of gains or losses \citep{krokhmal_portfolio_2002},
rather than choosing a utility function. In addition, CVaR is a coherent
risk measure and thus has strong decision-theoretic support \citep{artzner_coherent_1999}.
Furthermore, CVaR can be formulated as a convex optimization problem
and it is thus more tractable \citep{rockafellar_optimization_2000}
compared to other risk-aware objectives. \citet{maceira_application_2015}
applied CVaR to a multi-stage power generation planning problem, and
solved it by combining a scenario-tree based method with stochastic
dual dynamic programming. Applications of CVaR in capacity planning
problems can be found in \citep{szolgayova_dynamic_2011,delgado_transmission_2013}.

\subsection{Solution Methods}

Early work on MCEP modeled the problem as a Markov decision process
(MDP) and solved it with (exact) dynamic programming \citep{wu_innovative_2010,wu_efficient_2012}
or approximate dynamic programming \citep{zhao_approximate_2017}.
These methods are subject to the curse of dimensionality; more specifically,
the size of the action space of the MDPs grows exponentially in the
number of facilities. Alternative solution methods for risk-averse
MDP can be found in \citep{ruszczynski_risk-averse_2010,haskell_convex_2015},
but these methods can be inefficient when the actions are discrete
and high-dimensional.

In an alternative stream, scenario tree-based multi-stage stochastic
programming has been widely applied to MCEP \citep{huang_value_2009,taghavi_multi-stage_2016}.
In this method, the evolution of uncertain parameters is modeled as
a scenario-tree, and the model is solved with decomposing by fixing
the allocation plan \citep{huang_value_2009}, by a Benders decomposition-based
heuristic \citep{taghavi_multi-stage_2016}, or by Lagrangian relaxation
\citep{taghavi_lagrangian_2018}. Nevertheless, the size of the scenario
tree grows exponentially when the number of stages or the dimension
of uncertain parameters increases.

To address this open problem, decision rule-based multi-stage stochastic
programming was proposed. This method approximates the policy space
with parameterized decision rules. Then, the focus is on optimizing
the parameters of the decision rule rather than optimizing the policy
itself. Well-known decision rules include linear and piecewise linear
\citep{georghiou_generalized_2015}, but these rules may not be applicable
to MCEP as the system is usually modular---i.e., the capacity is
discrete. To solve MCEPs with discrete capacity, researchers have
investigated if--then decision rules and proposed customized decomposition
algorithms to optimize the parameters \citep{cardin_approach_2017,zhao_decision_2018}.
The solution technique in our present paper is similar to the branch-and-cut
based decomposition (BAC-D) algorithm proposed in \citep{zhao_decision_2018},
but our study differs from the previous literature in the following
two respects. First, in our paper, we construct cuts for the master
problem with the subgradients of the recourse function. These subgradients
can be computed by solving some small-scale linear programs (LPs).
Compared to the BAC-D algorithm that constructs cuts by solving large-scale
LPs, our method is much less time-consuming and it can also provide
more interpretation for the resulting policy. Second, both BAC-D and
the subgradient method cannot ensure convergence to the global optimum,
but we discuss how to improve the best-found solution by using a multi-cut
version of our algorithm.

\section{\label{sec:Model-Formulation}Model Formulation}

For our multi-stage capacity expansion problem, we introduce a set
of facilities $\mathcal{N}\triangleq\left\{ 1,\ldots,N\right\} $,
a set of customers $\mathcal{I}\triangleq\left\{ 1,\ldots,I\right\} $,
and a finite planning horizon $\mathcal{T}\triangleq\left\{ 1,\ldots,T\right\} $.
In each period, the demand generated by the customers is allocated
to facilities subject to available capacity. Given observations of
demand, the decision maker can expand the capacity at the end of each
period. The objective is to maximize the net present value (ENPV)
of the system by optimizing the capacity expansion policy. In this
paper, we focus on capacity planning at the strategic level, so the
following assumptions are made: 1) there is no temporary shutdown
of the facilities and the contraction of capacity is not allowed;
2) the expansion lead time is negligible. The notation for our model
is summarized in Table \ref{tab:Notations_MCEP}.

\begin{table}
\caption{\label{tab:Notations_MCEP}Notations for the multi-stage MCEP}

\centering{}{\small{}}%
\begin{tabular}{ll}
\hline 
{\small{}$\mathcal{I}$} & {\small{}Set of customers, $i\in\mathcal{I}$ and $\left|\mathcal{I}\right|=I$}\tabularnewline
{\small{}$\mathcal{N}$} & {\small{}Set of facilities, $n\in\mathcal{N}$, and $\left|\mathcal{N}\right|=N$}\tabularnewline
{\small{}$\mathcal{T}$} & {\small{}Set of time periods, $t\in\mathcal{T}$, and $\left|\mathcal{T}\right|=T$}\tabularnewline
$\mathcal{L}$ & {\small{}Set of line segments of the piecewise expansion cost, $l\in\mathcal{L}$
and $\left|\mathcal{L}\right|=L$}\tabularnewline
{\small{}$\Xi_{t}$} & {\small{}Sample space of the uncertain demands in time $t\in\mathcal{T}$}\tabularnewline
{\small{}$\xi_{it}$} & {\small{}Amount of demand generated from customer $i\in\mathcal{I}$
in time $t\in\mathcal{T}$; its vector form is}\tabularnewline
 & {\small{}$\xi_{t}=\left(\xi_{1t},\ldots,\xi_{It}\right)$ such that
$\xi_{t}\in\Xi_{t}$}\tabularnewline
{\small{}$K_{n}^{\max}$} & {\small{}The maximum capacity of facility $n\in\mathcal{N}$ that
can be installed; its vector form is}\tabularnewline
 & {\small{}$K^{\max}=\left(K_{1}^{\max},\ldots,K_{N}^{\max}\right)$}\tabularnewline
{\small{}$\mathbb{K}$} & {\small{}The feasible set of capacity $\mathbb{K}\triangleq\left\{ \left.K\in\mathbb{Z}_{+}^{N}\right|K\le K^{\max}\right\} $}\tabularnewline
{\small{}$K_{nt}$} & {\small{}Capacity of facility $n\in\mathcal{N}$ in time $t\in\mathcal{T}\cup\left\{ 0\right\} $;
its vector form is $K_{t}=\left(K_{1t},\ldots,K_{Nt}\right)$}\tabularnewline
{\small{}$\mathcal{K}_{t}$} & {\small{}Policies that from the historical demands $\xi_{\left[t\right]}$
to the capacity in time $t\in\mathcal{T}$}\tabularnewline
{\small{}$\mathcal{K}$} & {\small{}Set of policies that are feasible to MCEP}\tabularnewline
{\small{}$\gamma$} & {\small{}Discount factor, $0<\gamma<1$}\tabularnewline
{\small{}$r_{int}$} & {\small{}Unit revenue from satisfying customer $i\in\mathcal{I}$
with facility $n\in\mathcal{N}$ in time $t\in\mathcal{T}$}\tabularnewline
{\small{}$b_{it}$} & {\small{}Unit penalty cost for unsatisfied customer $i\in\mathcal{I}$
in time $t\in\mathcal{T}$}\tabularnewline
{\small{}$p_{nlt}$} & {\small{}Slope of the $l^{th}$ line segment of the expansion costs
corresponding to facility $n$ in time $t$}\tabularnewline
{\small{}$q_{nlt}$} & {\small{}Intercept of the $l^{th}$ line segment of the expansion
costs corresponding to facility $n$ in time $t$}\tabularnewline
{\small{}$z_{int}$} & {\small{}Amount of demands allocated from customer $i\in\mathcal{I}$
to facility $n\in\mathcal{N}$ in time $t\in\mathcal{T}$}\tabularnewline
\hline 
\end{tabular}{\small\par}
\end{table}

Let $\mathbb{K}\triangleq\left\{ K\in\mathbb{Z}_{+}^{N}|K\le K^{\max}\right\} $
be the finite set of possible capacity levels, where $K^{\max}\triangleq\left(K_{1}^{\max},\ldots,K_{N}^{\max}\right)$
is the vector of maximum possible capacity levels. Denote $K_{t}\triangleq\left(K_{1t},\ldots,K_{Nt}\right)\in\mathbb{K}$
as the vector of the installed capacity at the end of time period
$t\in\mathcal{T}\cup\left\{ 0\right\} $, and $K_{\left[t\right]}\triangleq\left(K_{0},K_{1}\ldots,K_{t}\right)$
as the history of installed capacity levels up to time $t$. Denote
$\Delta K_{nt}\triangleq K_{nt}-K_{n\left(t-1\right)}$ and $\Delta K_{t}\triangleq\left(\Delta K_{1t},\ldots,\Delta K_{Nt}\right)$
as the change in capacity at time $t\in\mathcal{T}$. Without loss
of generality, we assume that no capacity is installed at the beginning
so that $\Delta K_{0}=K_{0}$.

Denote $\Xi_{t}\subset\mathbb{R}^{I}$ as the sample space for the
customer demand in time $t\in\mathcal{T}\cup\left\{ 0\right\} $,
and $\xi_{t}\triangleq\left(\xi_{1t},\ldots,\xi_{It}\right)$ as the
realized demand such that $\xi_{t}\in\Xi_{t}$. Without loss of generality,
we assume that the demand at $t=0$, i.e. $\xi_{0}\in\Xi_{0}$, is
known. We further denote $\xi_{\left[t\right]}\triangleq\left(\xi_{0},\xi_{1},\ldots,\xi_{t}\right)$
as the history of demand up to time $t$, $\Xi\triangleq\times_{t=0}^{T}\Xi_{t}$
to be the set of all possible demand realizations, and $\xi\triangleq\xi_{\left[T\right]}$.

\subsection{Policies for Flexible MCEPs}

In the flexible MCEP, capacity decisions are made sequentially based
on observations of demand. Specifically, $K_{t}$ is a mapping from
historical demand $\xi_{\left[t\right]}$ to capacity expansion decisions,
for all $t\in\mathcal{T}$. We denote the policies for the flexible
MCEP as
\[
\mathcal{K}_{t}:\Xi_{0}\times\cdots\times\Xi_{t}\mapsto\mathbb{K},\ \ \forall t=0,\ldots,T.
\]
We require that $\mathcal{K}_{t}$ for all $t\in\mathcal{T}\cup\left\{ 0\right\} $
be \emph{non-anticipativ}e; the capacity decision $K_{t}=\mathcal{K}_{t}\left(\xi_{\left[t\right]}\right)$
in time $t$ may only depend on $\xi_{\left[t\right]}$ (it does not
have access to future information).

We further denote $\mathcal{K}_{\left[t\right]}\triangleq\left(\mathcal{K}_{0},\ldots,\mathcal{K}_{t}\right)$
as the policies up to time $t$. As contraction of capacity is not
allowed, the set of feasible policies for the flexible MCEP is
\[
\bar{\mathcal{K}}\triangleq\left\{ \left.\left(\mathcal{K}_{0},\ldots,\mathcal{K}_{T}\right)\right|\mathcal{K}_{t-1}\left(\xi_{\left[t-1\right]}\right)\le\mathcal{K}_{t}\left(\xi_{\left[t\right]}\right),\forall t\in\mathcal{T},\xi\in\Xi\right\} .
\]

\subsection{Profits and Costs}

Denote $\Pi_{t}\left(K_{t-1},\xi_{t}\right)$ as the profit given
the realized demand $\xi_{t}$ and the installed capacity $K_{t-1}$.
Denote $z_{int}$ as the demand allocated from customer $i\in\mathcal{I}$
to facility $n\in\mathcal{N}$ in time $t\in\mathcal{T}$. Then, $\Pi_{t}\left(K_{t-1},\xi_{t}\right)$
is given by the value of the following linear program:{\small{}
\begin{alignat}{2}
\Pi_{t}\left(K_{t-1},\xi_{t}\right)\triangleq\max_{z}\  & \sum_{i\in\mathcal{I}}\sum_{n\in\mathcal{N}}r_{int}z_{int}-\sum_{i\in\mathcal{I}}b_{it}\left(\xi_{it}-\sum_{n\in\mathcal{N}}z_{int}\right)\label{eq:Allocation_model}\\
\text{s.t.\  } & \sum_{n\in\mathcal{N}}z_{int}\le\xi_{it}, &  & \ \ \forall i\in\mathcal{I},\nonumber \\
 & \sum_{i\in\mathcal{I}}z_{int}\le K_{n\left(t-1\right)}, &  & \ \ \forall n\in\mathcal{N},\nonumber \\
 & z_{int}\ge0, &  & \ \ \forall i\in\mathcal{I},n\in\mathcal{N},\nonumber 
\end{alignat}
}where $r_{int}$ is the unit revenue for satisfying customer $i$'s
demand with facility $n$, and $b_{it}$ is the unit penalty for unsatisfied
demand of customer $i$. We do not enforce the constraint that all
demands must be met.

Denote $c_{t}\left(\Delta K_{t}\right)$ for all $t\in\mathcal{T}\cup\left\{ 0\right\} $
as the capacity expansion cost given $\Delta K_{t}$, we assume $c_{t}\left(\cdot\right)$
is piecewise linear. Denote $\mathcal{L}\triangleq\left\{ 1,\ldots,L\right\} $
as a set of indices for $L$ line segments, and denote $\left(a_{n1},\ldots,a_{n\left(L+1\right)}\right)$
as a set of breakpoints for the expansion costs for facility $n\in\mathcal{N}$
such that $a_{n1}=0$ and $K_{n}^{\max}<a_{n\left(L+1\right)}$. Let
$p_{nlt}$ and $q_{nlt}$ be the slope and intercept of the $l^{th}$
line segment of the expansion costs for facility $n$ in time $t$.
The cost function is then:{\small{}
\begin{multline}
c_{t}\left(\Delta K_{t}\right)\triangleq\left\{ \left.\sum_{n\in\mathcal{N}}c_{nt}\left(\Delta K_{nt}\right)\right|c_{nt}\left(\Delta K_{nt}\right)=p_{nlt}\cdot\Delta K_{nt}+q_{nlt},\ \text{if }\Delta K_{nt}\in\left[a_{nl},a_{n\left(l+1\right)}\right),\ \forall l\in\mathcal{L}\right\} ,\\
\ \forall t\in\mathcal{T}\cup\left\{ 0\right\} .\label{eq:linearize_cost_func}
\end{multline}
}Note that $c_{t}\left(\cdot\right)$ can be \emph{concave} in $\Delta K_{t}$,
as the expansion costs may enjoy the economies of scale. Eq. (\ref{eq:linearize_cost_func})
can represent/approximate a variety of concave cost functions; for
example, the fixed-charge function and the power function \citep{van_mieghem_commissioned_2003}.
\begin{example}
The fixed-charge function addresses economies of scale---that is,
the system incurs a fixed cost whenever capacity is expanded, and
the remaining costs increase linearly. We set the number of line segments
to be two and take $\left[a_{n1},a_{n2}\right)=\left[0,1\right)$
and $\left[a_{n2},a_{n3}\right)=\left[1,K_{n}^{\max}+1\right)$, and
$p_{n2}<p_{n1}$ and $q_{n2}>q_{n1}=0$.
\end{example}
\begin{example}
The power cost function can be represented by $\left(\Delta K_{nt}\right)^{v}$
where $0\le v<1$ is the factor for the economies of scale. It can
be approximated by piecewise linear functions; in addition, the relative
gap between the original cost function and the approximate one is
smaller than $1.5\%$ if we take the line segments and breakpoints
to be $\left(a_{1},a_{2},a_{3},\ldots,a_{L+1}\right)=\left(0,2^{0},2^{1},\ldots,2^{L-1}\right)$
for $L$ satisfying $2^{L-2}\le\max\left\{ K_{1}^{\max},\ldots,K_{N}^{\max}\right\} \le2^{L-1}$
\citep{zhao_decision_2018}.
\end{example}
\begin{rem}
Essentially, Eq. (\ref{eq:linearize_cost_func}) can formulate arbitrary
proper cost functions. As $\Delta K_{nt}$ is finite (i.e. $\Delta K_{nt}\in\left\{ 0,1,\ldots,K_{n}^{\text{max}}\right\} $),
we can set $L=K_{n}^{\max}+1$ so that each line segment corresponds
to a specific expansion cost at point $\Delta K_{nt}$.
\end{rem}
We denote the discount factor as $0<\gamma\le1$. Given policies $\mathcal{K}_{\left[T\right]}\in\mathcal{\bar{\mathcal{K}}}$
and the profit/cost structure described above, the cumulative future
costs from stages $t=0$ to $t=T$ for a particular $\xi\in\Xi$ are{\small{}
\[
Q\left(\mathcal{K}_{\left[T\right]},\xi\right)\triangleq\left\{ \left.c_{0}\left(K_{0}\right)+\sum_{t=1}^{T}\gamma^{t}\left(c_{t}\left(K_{t}-K_{t-1}\right)-\Pi_{t}\left(K_{t-1},\xi_{t}\right)\right)\right|K_{t}=\mathcal{K}_{t}\left(\xi_{\left[t\right]}\right),\forall t\in\mathcal{T}\cup\left\{ 0\right\} ,\xi\in\Xi\right\} .
\]
}The overall objective of the MCEP is to minimize $Q\left(\mathcal{K}_{\left[T\right]},\xi\right)$
(or equivalently, to maximize $-Q\left(\mathcal{K}_{\left[T\right]},\xi\right)$).

\subsection{A Risk-Averse Flexible MCEP}

In this subsection, we first formally present the risk-neutral flexible
MCEP and then its risk-averse counterpart. If the decision maker is
risk-neutral, the objective of the flexible MCEP is to find the expansion
policy that maximizes the ENPV:
\begin{alignat}{1}
\text{ENPV}_{\text{flex}}\triangleq\max_{\mathcal{K}_{\left[T\right]}\in\mathcal{\bar{\mathcal{K}}}} & \mathbb{E}\left[-Q\left(\mathcal{K}_{\left[T\right]},\xi\right)\right].\label{eq:Risk_Neutral_Model}
\end{alignat}

The costs in the above model can be especially high for some particular
realizations of $\xi$, especially when the variances of demand is
high. We incorporate CVaR into the objective as a remedy, the definition
of CVaR is next.

\renewcommand{\labelenumi}{(\roman{enumi})}
\begin{defn}
\citep{chen_value-at-risk_2008} Denote $X$ as a continuous random
variable and $F_{X}\left(y\right)=P\left\{ X\le y\right\} $ for all
$y\in\mathbb{R}$ as its cumulative distribution function.
\begin{enumerate}
\item The VaR of $X$ at confidence level $\alpha\in\left(0,1\right)$ is
\[
\text{VaR}_{\alpha}\left(X\right)\triangleq\inf\left\{ y|F_{X}\left(y\right)\ge\alpha\right\} .
\]
\item The CVaR of $X$ at confidence level $\alpha\in\left(0,1\right)$
is
\[
\text{CVaR}_{\alpha}\left(X\right)\triangleq\int_{-\infty}^{\infty}ydF_{X}^{\alpha}\left(y\right),
\]
where 
\[
F_{X}^{\alpha}\left(y\right)\triangleq\begin{cases}
0, & \text{when }y<\text{VaR}_{\alpha}\left(X\right),\\
\frac{F_{X}\left(y\right)-\alpha}{1-\alpha}, & \text{when }y\ge\text{VaR}_{\alpha}\left(X\right).
\end{cases}
\]
\end{enumerate}
\end{defn}
\renewcommand{\labelenumi}{\arabic{enumi}.} For a continuous random
variable $X$ and a confidence level $\alpha\in\left(0,1\right)$,
$\text{CVaR}_{\alpha}\left(X\right)$ is the conditional expectation
of $X$ greater than or equal to $\text{VaR}_{\alpha}\left(X\right)$.
This coincides the definition of ``expected shortfall'' \citep{acerbi_spectral_2002}.
If $\alpha\rightarrow1$, $\text{CVaR}_{\alpha}\left(X\right)$ approaches
the worst-case cost; whereas if $\alpha\rightarrow0$, $\text{CVaR}_{\alpha}\left(X\right)$
approaches the expectation of $X$.

We introduce a mean-CVaR objective for the MCEP by introducing a weight
factor $0\le\beta\le1$ \citep{shapiro_analysis_2011} to obtain:
\begin{equation}
\max_{K_{\left[T\right]}\in\bar{\mathcal{K}}}\ \beta\mathbb{E}\left[-Q\left(\mathcal{K}_{\left[T\right]},\xi\right)\right]+\left(1-\beta\right)\left[-\text{CVaR}_{\alpha}\left(Q\left(\mathcal{K}_{\left[T\right]},\xi\right)\right)\right].\label{eq:Obj_E_CVaR}
\end{equation}
In this formulation, decision makers can compromise between risk-neutral
and risk-averse policies by adjusting the weight factor $\beta$.
If we choose $\beta=1$, we recover the original risk-neutral model;
conversely, if $\beta=0$, we minimize $\text{CVaR}_{\alpha}\left(Q\left(\mathcal{K}_{\left[T\right]},\xi\right)\right)$.
\begin{rem}
\label{rem:spectral_risk_mesure}The objective in Eq. (\ref{eq:Obj_E_CVaR})
can be understood as a weighted sum of two ``expected shortfall''
measures (since the expectation is an expected shortfall for $\alpha=0$).
Expected shortfall is a spectral risk measure \citep{acerbi_spectral_2002},
so the objective of Eq. (\ref{eq:Obj_E_CVaR}) is essentially a special
case of the spectral risk measure (which is also coherent).
\end{rem}
Based on \citep[Theorem 2]{rockafellar_optimization_2000}, we can
introduce an auxiliary variable $u\in\mathbb{R}$ and calculate $\text{CVaR}_{\alpha}\left(Q\left(\mathcal{K}_{\left[T\right]},\xi\right)\right)$
by solving the following optimization problem: 
\begin{alignat*}{1}
\text{CVaR}_{\alpha}\left(Q\left(\mathcal{K}_{\left[T\right]},\xi\right)\right) & =\inf_{u\in\mathbb{R}}\left\{ u+\frac{1}{1-\alpha}\mathbb{E}\left[Q\left(\mathcal{K}_{\left[T\right]},\xi\right)-u\right]_{+}\right\} ,\\
 & =-\sup_{u\in\mathbb{R}}\left\{ -u-\frac{1}{1-\alpha}\mathbb{E}\left[Q\left(\mathcal{K}_{\left[T\right]},\xi\right)-u\right]_{+}\right\} ,
\end{alignat*}
where $\left[\cdot\right]_{+}$ denotes $\max\left\{ \cdot,0\right\} $.
Problem (\ref{eq:Obj_E_CVaR}) is then equivalent to
\begin{equation}
\text{ENPV}_{\alpha}\left(\beta\right)\triangleq\max_{u\in\mathbb{R},K_{\left[T\right]}\in\bar{\mathcal{K}}}\ -\left(1-\beta\right)u-\mathbb{E}\left[\beta Q\left(\mathcal{K}_{\left[T\right]},\xi\right)+\frac{1-\beta}{1-\alpha}\left[Q\left(\mathcal{K}_{\left[T\right]},\xi\right)-u\right]_{+}\right],\label{eq:Risk_aversion_Model}
\end{equation}
using the variational form of CVaR. Like Problem (\ref{eq:Risk_Neutral_Model}),
Problem (\ref{eq:Risk_aversion_Model}) is a multi-stage stochastic
programming problem. If $\beta=1$, Problem (\ref{eq:Risk_Neutral_Model})
is equivalent to Problem (\ref{eq:Risk_aversion_Model}). In addition,
we have the following result since $\text{CVaR}_{\alpha}\left(X\right)\ge\mathbb{E}\left[X\right]$
for all $\alpha\in\left(0,\,1\right)$ for any continuous random variable
$X$.

\renewcommand{\labelenumi}{(\roman{enumi})}
\begin{prop}
\label{prop:ENPV_propertise}Given Problems (\ref{eq:Risk_Neutral_Model})
and (\ref{eq:Risk_aversion_Model}):
\end{prop}
\begin{enumerate}
\item \emph{$\text{ENPV}_{\alpha}\left(\beta\right)$ is non-decreasing
in $\beta$ given any $\alpha\in\left(0,1\right)$.}
\item \emph{$\text{ENPV}_{\text{flex}}\ge\text{ENPV}_{\alpha}\left(\beta\right)$
for any $\beta\in\left[0,1\right],\alpha\in\left(0,1\right)$.}
\end{enumerate}
\renewcommand{\labelenumi}{\arabic{enumi}.}

Proposition \ref{prop:ENPV_propertise} states that $\text{ENPV}_{\alpha}\left(\beta\right)$
will not exceed $\text{ENPV}_{\text{flex}}$ given appropriately chosen
$\alpha$ and $\beta$. Furthermore, as $\beta$ decreases, decision
makers put more emphases on minimizing CVaR rather than minimizing
expected cost, and so $\text{ENPV}_{\alpha}\left(\beta\right)$ decreases.
Our interpretation is that, as the decision maker becomes more risk-averse,
it will tend to underestimate system performance.

From another perspective, we often want to estimate the value of flexibility
over an inflexible benchmark problem (i.e. the system has no capacity
adjustment options, detailed discussions on this problem appear in
Section \ref{sec:Numerical-Study}). Problems (\ref{eq:Risk_Neutral_Model})
and (\ref{eq:Risk_aversion_Model}) have the same level of flexibility
(both have capacity adjustment options and facility switching options),
and so the willingness of the decision maker to enable flexibility
may decrease as risk-aversion increases.

\section{\label{sec:Approximation}Decision Rules}

Problem (\ref{eq:Risk_aversion_Model}) is a multi-stage stochastic
programming problem with a non-convex objective. It is widely believed
that multi-stage stochastic programming is in general ``computationally
intractable already when medium-accuracy solutions are sought'' \citep{shapiro_complexity_2005}.

To develop a tractable solution strategy, we approximate $\mathcal{\overline{K}}$
with decision rules. That is, we restrict the policy space to a class
of parameterized functions $\mathcal{\tilde{K}}\left(\Theta\right)$,
where $\Theta\subset\mathbb{R}^{\text{dim}\left(\Theta\right)}$ is
some admissible set of the parameters. Then, we can optimize the parameters
$\theta\in\Theta$ which determine the decision rule $\mathcal{\tilde{K}}\left(\Theta\right)$,
instead of optimizing over all non-anticipative policies in $\mathcal{\bar{\mathcal{K}}}$.
Of course, $\mathcal{\tilde{K}}\left(\Theta\right)\subset\mathcal{\bar{\mathcal{K}}}$,
and so a decision rule may not be optimal for the original problem.
However, decision rules offer significant computational advantages
as well as managerial insight. In particular, they are far more accessible
to non-experts in practice.

\subsection{If--Then Decision Rules}

We focus on if--then decision rules to approximate the policy space
$\bar{\mathcal{K}}$. The motivation for our choice of if--then decision
rules is threefold. First, expansion decisions are binary by their
very nature---the capacity is either expanded or not. Also, the output
of the decision rule should be integral since capacity is discrete,
so a nonlinear decision rule is required. Second, some optimal if--then
policies for capacity expansion problems have been reported in the
literature. For example, the invest-stay put-disinvest policy from
\citep{eberly_multi-factor_1997} and the (s, S) policy from \citep{angelus_optimal_2000}
are both optimal for their respective problems. Therefore, it is reasonable
to speculate that an optimal if-then decision rule should at least
have good performance from the perspective of Problem (\ref{eq:Risk_aversion_Model}).
Third, if--then decision rules mimic the behavior of human beings
and are more intuitive in practical implementation \citep{cardin_approach_2017}.

An if--then decision rule in a single facility setting is stated
as: if the capacity gap of the facility exceeds a threshold, then
we expand capacity up to a certain level. However, in multi-facility
problems the dimensions of capacity levels and the demands may not
be equal, and so we need to transform one to make it comparable to
the other in order to calculate the capacity gap. We take $W\in\mathbb{R}^{I\times N}$
as a preset weight matrix such that $W^{\top}\xi\in\mathbb{R}^{N}$.
Then, we may compute weighted capacity gaps $\left\lfloor W^{\top}\xi\right\rceil -K_{t-1}$,
where $\left\lfloor \cdot\right\rceil $ denotes rounding to the nearest
integer. These weighted capacity gaps are the trigger condition for
our if--then rules.

The design of the weight matrix $W$ is case-specific. From a managerial
point of view, the entry of the weight matrix $W_{in}$ can be interpreted
as the profit coefficient of customer $i$ with respect to facility
$n$. In other words, if the per unit demand from customer $i$ is
more profitable for facility $n$, then we should tend to allocate
more demand from customer $i$ to facility $n$ and $W_{in}$ should
be larger.
\begin{example}
If the numbers of facilities and customers are the same (i.e. $I=N$)
then there is a bijective map from the facility/customer to its most
profitable customer/facility counterpart, and the weight matrix can
be chosen as an $N\times N$ identity matrix. In this case, the weighted
capacity gap of a facility is calculated in terms of subtracting the
current capacity from the most profitable demand.
\end{example}
\begin{example}
For more general cases, $W$ can be designed by calculating and weighting
the profit coefficients (i.e. $r_{int}+b_{it}$) of the allocation
model $\Pi_{t}\left(\cdot\right)$. We refer interested readers to
\citep{zhao_decision_2018} for further discussion of this choice.
\end{example}
Denote $\theta_{1,n}$ as the capacity adjustment parameter and $\theta_{2,nt}$
as the threshold parameter, and define the parameter vectors $\theta_{1}\triangleq\left(\theta_{1,1},\ldots,\theta_{1,N}\right)$
and $\theta_{2}\triangleq\left(\theta_{2,11},\ldots,\theta_{2,NT}\right)$.
The admissible sets for parameters $\theta_{1}$ and $\theta_{2}$
are
\[
\Theta_{1}\triangleq\left\{ \theta_{1}\in\mathbb{Z}^{N}|0\le\theta_{1}\le\theta_{1}^{\max}\right\} \ \ \text{and}\ \ \Theta_{2}\triangleq\left\{ \theta_{2}\in\mathbb{R}^{N\times T}|0\le\theta_{2}\le\theta_{2}^{\max}\right\} ,
\]
respectively, and we further define $\theta\triangleq\left(\theta_{1},\theta_{2}\right)$
and $\Theta\triangleq\Theta_{1}\times\Theta_{2}$ for succinctness.
Then, we may define our if--then decision rule as follows: for all
$n\in\mathcal{N},t\in\mathcal{T},\xi\in\Xi,K_{t-1}\in\mathbb{K}$,
\begin{equation}
\tilde{\mathcal{K}}_{nt}\left(K_{n\left(t-1\right)},\xi_{t};\theta\right)\triangleq\begin{cases}
\left\lfloor \sum\limits _{i\in\mathcal{I}}W_{in}\xi_{it}\right\rceil +\theta_{1,n}, & \begin{array}{c}
\text{if }\left\lfloor \sum\limits _{i\in\mathcal{I}}W_{in}\xi_{it}\right\rceil -K_{n\left(t-1\right)}\ge\theta_{2,nt}\ \text{and}\ \\
\left\lfloor \sum\limits _{i\in\mathcal{I}}W_{in}\xi_{it}\right\rceil +\theta_{1,n}\le K_{n}^{\max},
\end{array}\\
K_{n\left(t-1\right)}, & \text{otherwise}.
\end{cases}\label{eq:Decision_Rule}
\end{equation}
This decision rule states that if the weighted capacity gap of facility
$n$ (i.e. $\left\lfloor \sum_{i\in\mathcal{I}}W_{in}\xi_{it}\right\rceil -K_{n\left(t-1\right)}$)
exceeds the threshold $\theta_{2,nt}$ \emph{and} the expanded capacity
does not exceed the maximum capacity, then we expand the capacity
of facility $n$ up to level $\left\lfloor \sum_{i\in\mathcal{I}}W_{in}\xi_{it}\right\rceil +\theta_{1,n}$.
Otherwise, the capacity is unchanged. Note that the parameters $\theta_{2,nt}$
are continuous and allowed to vary with time (the purpose of which
is to enable a better approximation of the true optimal policy). In
addition, $\theta_{1,n}$ is integral so the decision rule automatically
yields integral expansion decisions. Now with Policy (\ref{eq:Decision_Rule}),
we want to optimize the scenario-independent parameters $\theta_{1}$
and $\theta_{2}$ with respect to the mean-CVaR.

\subsection{Risk-Averse MCEP with Decision Rule}

Denote $\mathcal{\tilde{K}}_{\left[T\right]}\left(\theta\right)\triangleq\left(\mathcal{K}_{0},\tilde{\mathcal{K}}_{1}\left(\cdot;\theta\right),\ldots,\tilde{\mathcal{K}}_{T}\left(\cdot;\theta\right)\right)$
for $\theta\in\Theta$ as the vector of parameterized if--then decision
rules encoded by Policy (\ref{eq:Decision_Rule}). Note that $\mathcal{K}_{0}:\Xi_{0}\mapsto\mathbb{K}$
is the policy for the initial capacity decision $K_{0}$, and is independent
of $\theta$. Then, given Policy (\ref{eq:Decision_Rule}), Problem
(\ref{eq:Risk_aversion_Model}) can be approximated with the following
model
\begin{alignat}{1}
\max_{\theta\in\Theta,u\in\mathbb{R}}\  & -\left(1-\beta\right)u-\mathbb{E}\left[\beta Q\left(\mathcal{\tilde{K}}_{\left[T\right]}\left(\theta\right),\xi\right)+\frac{1-\beta}{1-\alpha}\left[Q\left(\mathcal{\tilde{K}}_{\left[T\right]}\left(\theta\right),\xi\right)-u\right]_{+}\right],\label{eq:Risk_Model_DR}
\end{alignat}
which optimizes over policies of the form (\ref{eq:Decision_Rule}).
We now want to optimize the policy parameters $\theta$, rather than
the policy $\mathcal{K}\in\mathcal{\bar{K}}$, which is more tractable.

\subsection{MILP Transformations of the Decision Rule-based Model}

To transform Problem (\ref{eq:Risk_Model_DR}) into a computationally
workable form, we apply sample average approximation and linearization
to transform Problem (\ref{eq:Risk_Model_DR}) into an MILP.

Denote $\mathcal{S}\triangleq\left\{ 1,\ldots,S\right\} $ as a set
of scenario indices, and $\left\{ \xi^{1},\ldots,\xi^{S}\right\} $
as a set of sample paths generated via Monte-Carlo simulation of customer
demand, where $\xi^{s}\triangleq\left(\xi_{1}^{s},\ldots,\xi_{T}^{s}\right)$
for all $s\in\mathcal{S}$. Now, we introduce a set of binary auxiliary
variables $\delta_{t}^{s}\triangleq\left(\delta_{1t}^{s},\ldots,\delta_{Nt}^{s}\right)$
and transform the if--then decision rule $\tilde{\mathcal{K}}_{t}\left(K_{t-1}^{s},\xi_{t}^{s};\theta\right)$
into constraints. Let $\delta_{nt}^{s}$ be a binary variable such
that for all $t\in\mathcal{T},s\in\mathcal{S}$, the capacity of facility
$n$ is expanded to $\left\lfloor \sum_{i\in\mathcal{I}}W_{in}\xi_{it}^{s}\right\rceil +\theta_{1,n}$
if $\delta_{nt}^{s}=1$; otherwise, the capacity is unchanged. This
effect is achieved by the following Big-M constraints:
\begin{subequations}
\label{Big_M_constraints}
\begin{alignat}{2}
\left\lfloor W^{\top}\xi_{t}^{s}\right\rceil -K_{t-1}^{s}-\theta_{2,t}\ge & \left(M\mathbb{I}_{N\times N}\right)^{\top}\left(\delta_{t}^{s}-\boldsymbol{1}_{N}\right), & \ \ \forall t\in\mathcal{T},s\in\mathcal{S},\label{eq:Big_M_constraints_1}\\
\left\lfloor W^{\top}\xi_{t}^{s}\right\rceil -K_{t-1}^{s}-\theta_{2,t}< & \left(M\mathbb{I}_{N\times N}\right)^{\top}\delta_{t}^{s}, & \ \ \forall t\in\mathcal{T},s\in\mathcal{S},\label{eq:Big_M_constraints_2}\\
K_{t}^{s}\le & \left\lfloor W^{\top}\xi_{t}\right\rceil +\theta_{1}+\left(M\mathbb{I}_{N\times N}\right)^{\top}\left(\boldsymbol{1}_{N}-\delta_{t}^{s}\right), & \ \ \forall t\in\mathcal{T},s\in\mathcal{S},\\
K_{t}^{s}\ge & \left\lfloor W^{\top}\xi_{t}^{s}\right\rceil +\theta_{1}+\left(M\mathbb{I}_{N\times N}\right)^{\top}\left(\delta_{t}^{s}-\boldsymbol{1}_{N}\right), & \ \ \forall t\in\mathcal{T},s\in\mathcal{S},\\
K_{t}^{s}-K_{t-1}^{s}\le & \left(M\mathbb{I}_{N\times N}\right)^{\top}\delta_{t}^{s}, & \ \ \forall t\in\mathcal{T},s\in\mathcal{S},\\
K_{t}^{s}\le & K^{\max}, & \ \ \forall t\in\mathcal{T},s\in\mathcal{S},\\
\delta_{t}^{s}\in & \left\{ 0,1\right\} ^{N}, & \ \ \forall t\in\mathcal{T},s\in\mathcal{S},\label{eq:Big_M_constraints_7}
\end{alignat}
\end{subequations}
where $M\gg0$ is a large constant, $\mathbb{I}_{N\times N}$ is an
$N\times N$ identity matrix, and $\boldsymbol{1}_{N}$ denotes an
$N$-dimensional vector of ones. In addition, the right-hand sides
of Policy (\ref{eq:Decision_Rule}) are integral so that Constraints
\eqref{Big_M_constraints} map from continuous capacity levels to
discrete capacity expansion decisions.

Note that the nonlinear term $\left[Q\left(\mathcal{\tilde{K}}_{\left[T\right]}\left(\theta\right),\xi\right)-u\right]_{+}$
in the objective of Problem (\ref{eq:Risk_Model_DR}). To linearize
it, we examine its epigraph by introducing auxiliary variables $\eta^{s}\in\mathbb{R}$
for all $s\in\mathcal{S}$, such that $\eta^{s}\ge Q\left(\mathcal{\tilde{K}}_{\left[T\right]}\left(\theta\right),\xi^{s}\right)-u$
and $\eta^{s}\ge0$. Problem\textcolor{blue}{{} }(\ref{eq:Risk_Model_DR})
is subsequently approximated via the aforementioned procedure to obtain:
\begin{subequations}
\label{saa_model}{\small{}
\begin{multline}
\min_{K_{0},u,\theta,\eta^{s},\delta^{s},K_{\left[T\right]}^{s}}\ \beta c_{0}\left(K_{0}\right)+\left(1-\beta\right)u+\frac{1}{S}\sum_{s=1}^{S}\left[\beta\sum_{t=1}^{T}\gamma^{t}\left(c_{t}\left(K_{t}^{s}-K_{t-1}^{s}\right)-\Pi_{t}\left(K_{t-1}^{s},\xi_{t}^{s}\right)\right)+\frac{1-\beta}{1-\alpha}\eta^{s}\right]\label{eq:problem_with_decision_rule}
\end{multline}
\vspace{-0.7in}
}{\small\par}

{\small{}
\begin{alignat}{2}
\text{s.t. } & \eqref{eq:Big_M_constraints_1}-\eqref{eq:Big_M_constraints_7} &  & \ \ \forall t\in\mathcal{T},s\in\mathcal{S},\label{eq:Problem_with_rule_(1)}\\
 & \eta^{s}\ge c_{0}\left(K_{0}\right)+\sum_{t=1}^{T}\gamma^{t}\left(c_{t}\left(K_{t}^{s}-K_{t-1}^{s}\right)-\Pi_{t}\left(K_{t-1}^{s},\xi_{t}^{s}\right)\right)-u,\ \  &  & \ \ \forall s\in\mathcal{S},\label{eq:Problem_with_rule_(2)}\\
 & \eta^{s}\in\mathbb{R}_{+},\ \ K_{t}^{s}\in\mathbb{R}^{N}, &  & \ \ \forall t\in\mathcal{T},s\in\mathcal{S},\label{eq:Problem_with_rule_(3)}\\
 & K_{0}^{s}=K_{0}, &  & \ \ \forall s\in\mathcal{S},\label{eq:Problem_with_rule_K0s}\\
 & K_{0}\in\mathbb{K},\ \ u\in\mathbb{R},\ \ \theta\in\Theta.\label{eq:Problem_with_rule_(4)}
\end{alignat}
}{\small\par}
\end{subequations}

Problem \eqref{saa_model} now appears as an MILP. In this MILP, Constraint
(\ref{eq:Problem_with_rule_K0s}) appears to simplify notation.

\section{\label{sec:BAC-D_Algorithm}Subgradient-based Decomposition Algorithm}

In this section, we propose a decomposition algorithm for Problem
\eqref{saa_model}. We first reformulate Problem \eqref{saa_model}
as a ``two-stage'' stochastic programming problem. Then, making
use of this structure, we decompose the ``second-stage'' over scenarios,
compute the sub-gradients of the recourse function, and then update
the parameters of the decision rule. This decomposition algorithm
may not converge to the global optimum, so we later introduce a multi-cut
procedure to improve upon the best-found solution of this algorithm.

\subsection{Two-Stage Decomposition}

Problem \eqref{saa_model} can be treated as a ``two-stage'' stochastic
programming problem. The first-stage decisions consist of the parameters
$\theta,$ the initial capacity $K_{0}$, and the auxiliary variable
$u$ from the variational form of CVaR. These decisions are all scenario-independent
(i.e. here-and-now decisions):
\begin{alignat}{1}
\min_{K_{0},u,\theta} & \left\{ \beta c_{0}\left(K_{0}\right)+\left(1-\beta\right)u+\frac{1}{S}\sum_{s=1}^{S}R^{s}\left(K_{0},u,\theta\right),\ \text{s.t. }\eqref{eq:Problem_with_rule_(4)}\right\} ,\label{eq:first_stage_problem}
\end{alignat}
where $R^{s}\left(K_{0},u,\theta\right)$ is the recourse function
of the second-stage problem on scenario $s\in\mathcal{S}$. The recourse
function here actually returns the multi-period revenue over time
periods $t=1,\ldots,T$. It determines the capacity plan $\left(K_{1}^{s},\ldots,K_{T}^{s}\right)$
and the auxiliary variables $\eta^{s}$, which are scenario-dependent
(i.e. wait-and-see decisions):{\small{}
\begin{alignat}{1}
R^{s}\left(K_{0},u,\theta\right)\triangleq\min_{\eta^{s},\delta^{s},K_{\left[T\right]}^{s}} & \left\{ \beta\sum_{t=1}^{T}\gamma^{t}\left(c_{t}\left(K_{t}^{s}-K_{t-1}^{s}\right)-\Pi_{t}\left(K_{t-1}^{s},\xi_{t}^{s}\right)\right)+\frac{1-\beta}{1-\alpha}\eta^{s},\ \text{s.t. }\eqref{eq:Problem_with_rule_(1)}\text{\textendash}\eqref{eq:Problem_with_rule_K0s}\right\} .\label{eq:second_stage_problem}
\end{alignat}
}This two-stage decomposition, however, is difficult to solve because
the recourse function $R^{s}\left(K_{0},u,\theta\right)$ is nonconvex
due to the integer wait-and-see variables (i.e. $\delta_{t}^{s}$)
and the nonconvex costs $c_{t}\left(\cdot\right)$.

Fortunately, this problem has appealing structure. Once the initial
capacity $K_{0}$, the parameters $\theta$ of Policy (\ref{eq:Decision_Rule}),
and $u$ are all fixed, we can determine $K_{\left[T\right]}^{s}\triangleq\left(K_{0},K_{1}^{s},\ldots,K_{T}^{s}\right)$
via Policy (\ref{eq:Decision_Rule}) given the demand vector $\xi^{s}$
on scenario $s\in\mathcal{S}$. Once $K_{\left[T\right]}^{s}$ is
known, the nonconvex costs $c_{t}\left(\cdot\right)$ for all $t\in\mathcal{T}$
can all be determined by Eq. (\ref{eq:linearize_cost_func}). Therefore,
Problem (\ref{eq:second_stage_problem}) can be solved, and the subgradients
of the recourse functions $R^{s}\left(K_{0},u,\theta\right)$ with
respect to $\left(K_{0},u,\theta_{1}\right)$ can be computed.

We remark that the subgradients of $R^{s}\left(K_{0},u,\theta\right)$
with respect to $\theta_{2}$ may not be attainable because of the
trigger conditions in Policy (\ref{eq:Decision_Rule}). For counter-example,
consider $\theta_{2,nt}=10$ and a simple if--then rule: if $\xi_{it}\ge\theta_{2,nt}$,
then expand the capacity of facility $n$ to $10$. In this case,
the dependence of the recourse function on $\theta_{2,nt}$ is implicit
since the rule is triggered for any $\xi_{it}>10$, and the subgradient
with respect to $\theta_{2,nt}$ is not available.

To address this difficulty, we take a primal decomposition and update
$\left(K_{0},u,\theta_{1}\right)$ and $\theta_{2}$ separately. To
update $\left(K_{0},u,\theta_{1}\right)$, we approximate the recourse
functions by using cut generation. Specifically, we solve the epigraph
formulation of Problem (\ref{eq:first_stage_problem}) in each iteration
$m$:

\vspace{-0.3in}

\begin{subequations}
\label{master_problem}
\begin{alignat}{2}
\min_{K_{0},u,\theta_{1},y}\  & \beta c_{0}\left(K_{0}\right)+\left(1-\beta\right)u+y\label{eq:problem_in_step3}\\
\text{s.t. } & y\ge\left(\phi^{m^{\prime}}\right)^{\top}\left(K_{0},u,\theta_{1}\right)+\phi_{0}^{m^{\prime}}, &  & \forall m^{\prime}=1,\ldots,m-1,\label{eq:Step_3_optimality_cut}\\
 & K_{0}\in\mathbb{K},\ \ u\in\mathbb{R},\ \ \theta_{1}\in\Theta_{1},\ \ y\in\mathbb{R},\label{eq:end_problem_in_step3}
\end{alignat}
\end{subequations}
where Eq. (\ref{eq:Step_3_optimality_cut}) contains the cuts generated
up to iteration $m-1$. Essentially, we try to use Eq. (\ref{eq:Step_3_optimality_cut})
to approximate the recourse functions from below. Problem \eqref{master_problem}
is an MILP with $N+N\times L$ integer variables, which can be directly
solved via commercial solvers.

\begin{figure}
\begin{centering}
\includegraphics[scale=0.9]{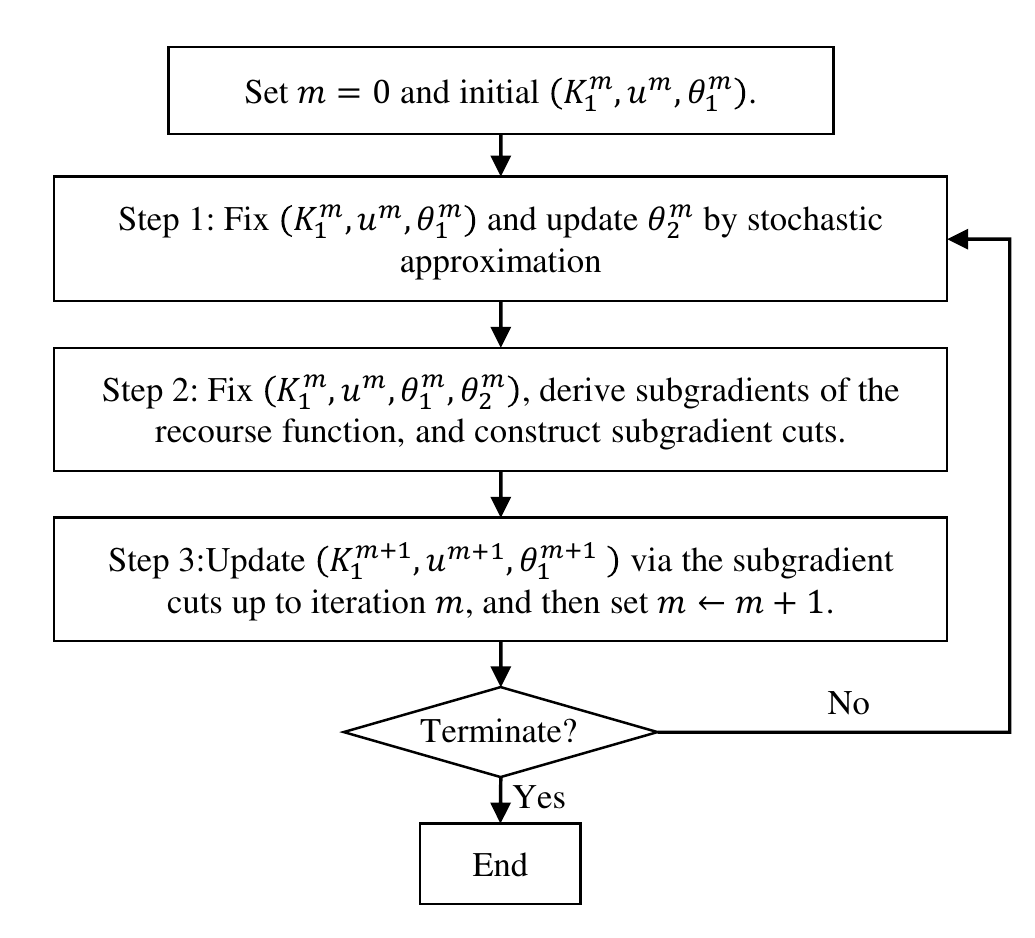}
\par\end{centering}
\caption{\label{fig:Algorithm_Flow_chart}Flow chart of the subgradient-based
decomposition algorithm.}
\end{figure}

Our algorithm contains three major steps. Denote $\left(K_{0}^{m},u^{m},\theta_{1}^{m}\right)$
as the optimal solution of Problem \eqref{master_problem}.
\begin{itemize}
\item Step 1: Fix $\left(K_{0}^{m},u^{m},\theta_{1}^{m}\right)$ and compute
$\theta_{2}^{m}$ by stochastic approximation. We construct and solve
single-scenario problems, and then iteratively average the resulting
optimal $\theta_{2}$.
\item Step 2: Once $\left(K_{0}^{m},u^{m},\theta_{1}^{m},\theta_{2}^{m}\right)$
are fixed, we compute subgradients of the recourse functions $R^{s}\left(\cdot\right)$
for all $s\in\mathcal{S}$ at $\left(K_{0}^{m},u^{m},\theta_{1}^{m}\right)$,
and construct subgradient cuts.
\item Step 3: We add the subgradient cut to Eq. (\ref{eq:Step_3_optimality_cut}),
and compute $\left(K_{0}^{m+1},u^{m+1},\theta_{1}^{m+1}\right)$ by
solving Problem \eqref{master_problem}. The algorithm then proceeds
iteratively until our termination conditions are met.
\end{itemize}
The overall framework of our algorithm is somewhat similar to the
one in \citep{zhao_decision_2018}. One of the main differences is
that in the second step, we solve small-scale LPs and compute the
subgradients analytically. This way is more efficient than solving
large-scale LPs. Figure \ref{fig:Algorithm_Flow_chart} gives a flow
chart of our algorithm. 

\subsection{Subgradient-based Decomposition Algorithm}

\subsubsection*{Step 1: Update $\text{\ensuremath{\theta_{2}}}$ via Stochastic Approximation}

Suppose we have fixed the first-stage decisions $\left(K_{0}^{m},u^{m},\theta_{1}^{m}\right)$
in the $m^{th}$ iteration, and let $k\ge0$ be a counter for the
iterations of stochastic approximation. We assign equal probability
to each scenario in $\mathcal{S}$, and randomly select one scenario
$\left\{ \xi_{1}^{s_{k}},\ldots,\xi_{T}^{s_{k}}\right\} $ for $s_{k}\in\mathcal{S}$
without replacement. Then, we construct the following single-scenario
problem given $\left(K_{0}^{m},u^{m},\theta_{1}^{m}\right)$ and $\xi^{s_{k}}$:
\begin{alignat}{2}
\min_{K_{\left[T\right]},\eta,\theta_{2}}\  & \beta\sum_{t=1}^{T}\gamma^{t}\left(c_{t}\left(K_{t}-K_{t-1}\right)-\Pi_{t}\left(K_{t-1},\xi_{t}^{s_{k}}\right)\right)+\frac{1-\beta}{1-\alpha}\eta+c_{\theta}\sum_{t\in\mathcal{T}}\sum_{n\in\mathcal{N}}\theta_{2,nt}\label{eq:ssp_model}\\
\text{s.t. } & K_{t}=\tilde{\mathcal{K}}_{t}\left(K_{t-1},\xi_{t}^{s_{k}};\left(\theta_{1}^{m},\theta_{2}\right)\right), &  & \ \ \forall t\in\mathcal{T},\nonumber \\
 & \eta\ge c_{0}\left(K_{0}^{m}\right)+\sum_{t=1}^{T}\gamma^{t}\left(c_{t}\left(K_{t}-K_{t-1}\right)-\Pi_{t}\left(K_{t-1},\xi_{t}^{s_{k}}\right)\right)-u^{m},\nonumber \\
 & \eta\in\mathbb{R}_{+},\ \ \theta_{2}\in\Theta_{2},\ \ K_{0}=K_{0}^{m},\ \ K_{\left[T\right]}\in\mathbb{R}^{N\times\left(T+1\right)},\nonumber 
\end{alignat}
where $c_{\theta}$ is a small positive constant intended to regularize
$\theta_{2}$ (this way the optimal $\theta_{2,nt}$ in Problem (\ref{eq:ssp_model})
is not a large number when facility $n$ is not expanded in time $t$).
The objective of the above problem is to find the optimal $\theta_{2}$
such that the cumulative costs given sample path $\left\{ \xi_{1}^{s_{k}},\ldots,\xi_{T}^{s_{k}}\right\} $
are minimized. Since $\left(K_{0}^{m},u^{m},\theta_{1}^{m}\right)$
are fixed and Policy (\ref{eq:Decision_Rule}) yields integral decisions,
the $K_{\left[T\right]}$ in Problem (\ref{eq:ssp_model}) are continuous
variables.
\begin{rem}
Note that in Problem (\ref{eq:ssp_model}), we have binary auxiliary
variables in $c_{t}\left(\cdot\right)$ and we need to transform $\tilde{\mathcal{K}}_{t}\left(\cdot\right)$
into Big-M constraints, so that Problem (\ref{eq:ssp_model}) is an
MILP with $N\times\left(L+1\right)\times T$ binary variables. For
modestly-sized problem instances, this MILP can be directly solved
by commercial solvers. For example, in our upcoming numerical study,
our problem has $N=5$, $L=5$, and $T=15$, and 450 binary variables.
In this case, Problem (\ref{eq:ssp_model}) can be solved by CPLEX
within seconds.
\end{rem}
Let $\hat{\theta}_{2}^{m,k}\triangleq\left(\hat{\theta}_{2,11}^{m,k},\ldots,\hat{\theta}_{2,NT}^{m,k}\right)$
be the optimal $\theta_{2}$ of Problem (\ref{eq:ssp_model}) given
sample path $\left\{ \xi_{1}^{s_{k}},\ldots,\xi_{T}^{s_{k}}\right\} $.
Then, we update $\theta_{2}$ by the following rule:
\begin{equation}
\bar{\theta}_{2}^{m,k+1}=\bar{\theta}_{2}^{m,k}+\sigma_{k}\left(\hat{\theta}_{2}^{m,k}-\bar{\theta}_{2}^{m,k}\right),\ \ \forall k\ge0,\label{eq:ssa_update_rule}
\end{equation}
where $\left\{ \sigma_{k}\right\} _{k\geq0}$ is a sequence of learning
rates that satisfies $\sum_{k=1}^{\infty}\sigma_{k}=\infty$ and $\sum_{k=1}^{\infty}\sigma_{k}^{2}=\infty$.
Given a preset precision and a minimum number of iterations $\underline{k}$,
we terminate the stochastic approximation update when $\left\Vert \bar{\theta}_{2}^{m,k+1}-\bar{\theta}_{2}^{m,k}\right\Vert _{\infty}\le\epsilon$
or $k\ge\bar{k}$, where $\left\Vert \cdot\right\Vert _{\infty}$
denotes the sup-norm. The above update rule is essentially the stochastic
approximation algorithm of \citep{robbins_stochastic_1951}. The advantage
is that we can derive an approximate $\theta_{2}$ by evaluating a
portion of, rather than all of, the scenarios in $\mathcal{S}$.

The update rule by Eq. (\ref{eq:ssa_update_rule}) may underestimate
$\theta_{2}$. Suppose the stochastic approximation algorithm terminates
after $k^{*}$ iterations. Denote $\hat{K}_{nt}^{m,k}$ for all $n\in\mathcal{N},t\in\mathcal{T}$
to be the optimal capacity decisions for Problem (\ref{eq:ssp_model})
given the selected scenario in the $k^{th}$ iteration, and let $\hat{\delta}_{nt}^{m,k}=1$
if $\hat{K}_{nt}^{m,k}>\hat{K}_{n\left(t-1\right)}^{m,k}$ and $\hat{\delta}_{nt}^{m,k}=0$
otherwise. Then, if for all of the selected scenarios, facility $n$
is not expanded in time $t$ (i.e. $\hat{\delta}_{nt}^{m,k^{\prime}}=0$
for all $k^{\prime}=1,\ldots,k^{*}-1$), the approximate $\bar{\theta}_{2,nt}^{m,k^{*}}$
may be an underestimate. This is because $\bar{\theta}_{2,nt}^{m,k^{*}}\le\max\left\{ \hat{\theta}_{2,nt}^{m,1},\ldots,\hat{\theta}_{2,nt}^{m,k^{*}-1}\right\} $
if it is updated via Eq. (\ref{eq:ssa_update_rule}). To fix this
issue, we take $M_{\theta}\gg0$ to be a large number such that $M_{\theta}\ge\hat{\theta}_{2,nt}^{m,k^{\prime}}$
for all $k^{\prime}=1,\ldots,k^{*}-1$. Then, the optimal $\theta_{2}$
given $\left(K_{0}^{m},u^{m},\theta_{1}^{m}\right)$ can be computed
by
\[
\theta_{2,nt}^{m}=\begin{cases}
M_{\theta} & \text{if }\sum\limits _{k^{\prime}=1}^{k^{*}-1}\hat{\delta}_{nt}^{m,k^{\prime}}=0,\\
\bar{\theta}_{2,nt}^{m,k^{*}}, & \text{otherwise,}
\end{cases}\ \ \forall n\in\mathcal{N},t\in\mathcal{T}.
\]

\subsubsection*{Step 2: Calculating the Subgradients of the Recourse Function}

Once $\theta_{2}^{m}$ is computed, we can fix $\left(K_{0}^{m},u^{m},\theta_{1}^{m},\theta_{2}^{m}\right)$
and determine the integer wait-and-see variables in Problem (\ref{eq:second_stage_problem}),
and then calculate the subgradients of the recourse functions with
respect to $\left(K_{0}^{m},u^{m},\theta_{1}^{m}\right)$. According
to Eq. (\ref{eq:second_stage_problem}), we need to compute the subgradients
of $c_{t}\left(\cdot\right)$ and $\Pi_{t}\left(\cdot\right)$ in
order to compute the subgradient of $R^{s}\left(K_{0}^{m},u^{m},\theta^{m}\right)$.
Therefore, in this step, we first derive the closed-form of the capacity
decisions $K_{t}$ with respect to $\left(K_{0}^{m},u^{m},\theta_{1}^{m}\right)$.
Subsequently, we calculate the subgradients $\partial c_{t}\left(\cdot\right)$
and $\partial\Pi_{t}\left(\cdot\right)$.

First, recall that $\delta_{nt}^{s}$ is the binary variable in Eq.
\eqref{Big_M_constraints} which indicates whether the capacity decision
of facility $n$ is triggered in time $t$. Given $\left(K_{0}^{m},u^{m},\theta_{1}^{m}\right)$
in the $m^{th}$ iteration and scenario $\left\{ \xi_{1}^{s},\ldots,\xi_{T}^{s}\right\} $
for $s\in\mathcal{S}$, the expansion decisions $\delta_{nt}^{m,s}$
and $K_{nt}^{m,s}$ can be determined sequentially from $t=1$ to
$T$ by:
\begin{equation}
\delta_{nt}^{m,s}=\begin{cases}
1, & \text{if }\left\lfloor \sum\limits _{i\in\mathcal{I}}W_{in}\xi_{it}^{s}\right\rceil -K_{n\left(t-1\right)}^{m,s}\ge\theta_{2,nt}^{m}\ \text{and}\ \left\lfloor \sum\limits _{i\in\mathcal{I}}W_{in}\xi_{it}^{s}\right\rceil +\theta_{1,n}^{m}\le K_{n}^{\max},\\
0, & \text{otherwise},
\end{cases}\ \ \forall t\in\mathcal{T},\label{eq:determine_expansion_decisions}
\end{equation}
and 
\[
K_{nt}^{m,s}=\delta_{nt}^{m,s}\left(\left\lfloor \sum\limits _{i\in\mathcal{I}}W_{in}\xi_{it}^{s}\right\rceil +\theta_{1}^{m}\right)+\left(1-\delta_{nt}^{m,s}\right)K_{n\left(t-1\right)}^{m,s},\ \ \forall t\in\mathcal{T},
\]
where we take $K_{n0}^{m,s}=K_{0}^{m}$ to simplify notation. We further
denote 
\begin{alignat*}{2}
h_{K,nt}^{m,s} & \triangleq\prod_{t_{0}=1}^{t}\left(1-\delta_{nt_{0}}^{m,s}\right), & \forall n\in\mathcal{N},t\in\mathcal{T}\cup\left\{ 0\right\} ,s\in\mathcal{S},\\
h_{\theta,nt}^{m,s} & \triangleq\sum_{t_{0}=1}^{t}\left[\prod_{t_{1}=t_{0}+1}^{t}\left(1-\delta_{nt_{1}}^{m,s}\right)\right]\delta_{nt_{0}}^{m,s}, & \forall n\in\mathcal{N},t\in\mathcal{T}\cup\left\{ 0\right\} ,s\in\mathcal{S},\\
h_{0,nt}^{m,s} & \triangleq\sum_{t_{0}=1}^{t}\left[\prod_{t_{1}=t_{0}+1}^{t}\left(1-\delta_{nt_{1}}^{m,s}\right)\right]\delta_{nt_{0}}^{m,s}\left\lfloor \sum\limits _{i\in\mathcal{I}}W_{in}\xi_{it_{0}}^{s}\right\rceil , & \ \ \forall n\in\mathcal{N},t\in\mathcal{T}\cup\left\{ 0\right\} ,s\in\mathcal{S}.
\end{alignat*}
For simplicity, we take $\prod_{t=t_{0}}^{t_{1}}x_{t}=1$ and $\sum_{t=t_{0}}^{t_{1}}x_{t}=0$
for any $x_{t}$ if $t_{0}>t_{1}$. Therefore, for $t=0$, we have
$h_{K,n0}^{m,s}=1$ and $h_{\theta,nt}^{m,s}=h_{0,n0}^{m,s}=0$. The
closed form of $K_{nt}^{m,s}$ with respect to $\left(K_{0}^{m},u^{m},\theta_{1}^{m}\right)$
follows.
\begin{lem}
\label{lem:closed_form_for_K}Given $\left(K_{0}^{m},u^{m},\theta^{m}\right)$,
scenario $s\in\mathcal{S}$, and Eq. (\ref{eq:determine_expansion_decisions}),
we have
\begin{equation}
K_{nt}^{m,s}=h_{K,nt}^{m,s}K_{n0}^{m}+h_{\theta,nt}^{m,s}\theta_{1,n}^{m}+h_{0,nt}^{m,s},\ \ \forall n\in\mathcal{N},t\in\mathcal{T}\cup\left\{ 0\right\} .\label{eq:closed_form_for_K}
\end{equation}
\end{lem}
Next, we compute the subgradients of the expansion costs. Note that
the expansion costs $c_{t}\left(\Delta K_{t}^{m,s}\right)$ for all
$t\in\mathcal{T}$ are deterministic since $K_{nt}^{m,s}$ are known.
Let $\Delta K_{nt}^{m,s}=K_{nt}^{m,s}-K_{n\left(t-1\right)}^{m,s}$
for all $t\in\mathcal{T}$ and $\Delta K_{n0}^{m,s}=K_{n0}^{m}$ for
all $s\in\mathcal{S}$, and define
\[
\tau_{nlt}^{m,s}\triangleq\begin{cases}
1, & \text{if }\Delta K_{nt}^{m,s}\in\left[a_{nl},a_{n\left(l+1\right)}\right),\\
0, & \text{otherwise},
\end{cases}\ \ \forall n\in\mathcal{N},l\in\mathcal{L},t\in\mathcal{T}\cup\left\{ 0\right\} .
\]
In the above, $\tau_{nlt}^{m,s}=1$ implies that $\Delta K_{nt}^{m,s}$
lies in the $l^{th}$ line segment of the expansion cost in Eq. (\ref{eq:linearize_cost_func}).
The slope and the intercept of the expansion cost for $\Delta K_{nt}^{m,s}$
are
\begin{align*}
g_{nt}^{m,s} & \triangleq\sum_{l=1}^{L}p_{nlt}\tau_{nlt}^{m,s},\ \ \text{and}\ \ g_{0,nt}^{m,s}\triangleq\sum_{l=1}^{L}q_{nlt}\tau_{nlt}^{m,s}, & \ \ \forall n\in\mathcal{N},t\in\mathcal{T}\cup\left\{ 0\right\} .
\end{align*}
Then, the subgradients of $c_{t}\left(\cdot\right)$ are as follows.
\begin{lem}
\label{lem:grad_C}Given $\left(K_{0}^{m},u^{m},\theta^{m}\right)$,
scenario $s\in\mathcal{S}$, and Eq. (\ref{eq:closed_form_for_K}),
a subgradient of $c_{t}\left(\cdot\right)$ for all $t\in\mathcal{T}\cup\left\{ 0\right\} $
at $\left(K_{0}^{m},u^{m},\theta_{1}^{m}\right)$ is 
\[
\left(\partial c_{t}^{m,s}\right)^{\top}\left(K_{0}^{m},u^{m},\theta_{1}^{m}\right)+\sum_{n\in\mathcal{N}}\left(g_{nt}^{m,s}h_{0,nt}^{m,s}-g_{nt}^{m,s}h_{0,n\left(t-1\right)}^{m,s}+g_{0,nt}^{m,s}\right),
\]
where
\[
\partial c_{t}^{m,s}=\left(\begin{array}{c}
\left(g_{nt}^{m,s}h_{K,nt}^{m,s}-g_{nt}^{m,s}h_{K,n\left(t-1\right)}^{m,s}\right)_{n\in\mathcal{N}}\\
0\\
\left(g_{nt}^{m,s}h_{\theta,nt}^{m,s}-g_{nt}^{m,s}h_{\theta,n\left(t-1\right)}^{m,s}\right)_{n\in\mathcal{N}}
\end{array}\right),\ \forall t\in\mathcal{T},\ \ \text{and}\ \ \partial c_{0}^{m,s}=\left(\begin{array}{c}
\left(g_{n0}^{m,s}h_{K,n0}^{m,s}\right)_{n\in\mathcal{N}}\\
0\\
0
\end{array}\right).
\]
\end{lem}
We now compute the subgradients of the profit. As $\Pi_{t}\left(K_{t-1},\xi_{t}\right)$
is given by an LP, the subgradient of $\Pi_{t}\left(K_{t-1},\xi_{t}\right)$
with respect to $\left(K_{0}^{m},u^{m},\theta_{1}^{m}\right)$ can
be computed from the dual to Problem (\ref{eq:Allocation_model}).
Let $\left(\mu_{nt}^{m,s}\right)_{n\in\mathcal{N},t\in\mathcal{T}}$
and $\left(\psi_{it}^{m,s}\right)_{i\in\mathcal{I},t\in\mathcal{T}}$
be the optimal dual variables with respect to the capacity and demand
constraints of $\Pi_{t}\left(\cdot\right)$, respectively.
\begin{lem}
\label{lem:grad_Pi}Given $\left(K_{0}^{m},u^{m},\theta^{m}\right)$,
scenario $s\in\mathcal{S}$, and Eq. (\ref{eq:closed_form_for_K}),
a subgradient of $\Pi_{t}\left(K_{n\left(t-1\right)}^{m,s},\xi_{it}^{s}\right)$
for all $t\in\mathcal{T}$ at $\left(K_{0}^{m},u^{m},\theta_{1}^{m}\right)$
is
\[
\left(\partial\Pi_{t}^{m,s}\right)^{\top}\left(K_{0}^{m},u^{m},\theta_{1}^{m}\right)+\sum_{i\in\mathcal{I}}b_{it}\xi_{it}^{s}-\sum_{i\in\mathcal{I}}\psi_{it}^{m,s}\xi_{it}^{s},
\]
where
\[
\partial\Pi_{t}^{m,s}=\left(\begin{array}{c}
\left(\mu_{nt}^{m,s}h_{K,n\left(t-1\right)}^{m,s}\right)_{n\in\mathcal{N}}\\
0\\
\left(\mu_{nt}^{m,s}h_{\theta,n\left(t-1\right)}^{m,s}\right)_{n\in\mathcal{N}}
\end{array}\right),\ \ \forall t\in\mathcal{T}.
\]
\end{lem}
Given Lemmas \ref{lem:closed_form_for_K}--\ref{lem:grad_Pi}, the
subgradients of the recourse functions can be computed via Proposition
\ref{prop:subgradient}.

\renewcommand{\labelenumi}{(\roman{enumi})}
\begin{prop}
\label{prop:subgradient}Given $\left(K_{0}^{m},u^{m},\theta^{m}\right)$
and scenario $s\in\mathcal{S}$, a subgradient in $\partial R^{s}\left(\cdot\right)$
with respect to $\left(K_{0}^{m},u^{m},\theta^{m}\right)$ is given
by:
\begin{enumerate}
\item If $c_{0}\left(K_{0}^{m}\right)+\sum_{t=1}^{T}\gamma^{t}\left(c_{t}\left(\Delta K_{t}^{s}\right)-\Pi_{t}\left(K_{t-1}^{s},\xi_{t}^{s}\right)\right)-u^{m}\ge0$,
we have
\begin{alignat*}{2}
\partial_{K_{0n}^{m}}R^{m,s} & =\beta\sum_{t=1}^{T}\gamma^{t}\left(g_{nt}^{m,s}\left(h_{K,nt}^{m,s}-h_{K,n\left(t-1\right)}^{m,s}\right)-\mu_{nt}^{m,s}h_{K,n\left(t-1\right)}^{m,s}\right), & \ \ \forall n\in\mathcal{N},\\
\partial_{u^{m}}R^{m,s} & =0,\\
\partial_{\theta_{1,n}^{m}}R^{m,s} & =\beta\sum_{t=1}^{T}\gamma^{t}\left(g_{nt}^{m,s}\left(h_{\theta,nt}^{m,s}-h_{\theta,n\left(t-1\right)}^{m,s}\right)-\mu_{nt}^{m,s}h_{\theta,n\left(t-1\right)}^{m,s}\right), & \forall n\in\mathcal{N},
\end{alignat*}
and
\[
R_{0}^{m,s}=\beta\sum_{t=1}^{T}\gamma^{t}\left(g_{nt}^{m,s}\left(h_{0,nt}^{m,s}-h_{0,n\left(t-1\right)}^{m,s}\right)+g_{0,t}^{m,s}+\sum_{i\in\mathcal{I}}b_{it}\xi_{it}^{s}-\sum_{i\in\mathcal{I}}\psi_{it}^{m,s}\xi_{it}^{s}-\mu_{nt}^{m,s}h_{0,nt}^{m,s}\right).
\]
\item If $c_{0}\left(K_{0}^{m}\right)+\sum_{t=1}^{T}\gamma^{t}\left(c_{t}\left(\Delta K_{t}^{s}\right)-\Pi_{t}\left(K_{t-1}^{s},\xi_{t}^{s}\right)\right)-u^{m}<0,$
we have{\small{}
\begin{alignat*}{2}
\partial_{K_{0n}^{m}}R^{m,s} & =\left(\beta+\frac{1-\beta}{1-\alpha}\right)\sum_{t=1}^{T}\gamma^{t}\left(g_{nt}^{m,s}\left(h_{K,nt}^{m,s}-h_{K,n\left(t-1\right)}^{m,s}\right)-\mu_{nt}^{m,s}h_{K,n\left(t-1\right)}^{m,s}\right)+\frac{1-\beta}{1-\alpha}g_{n0}^{m,s}, & \ \ \forall n\in\mathcal{N},\\
\partial_{u^{m}}R^{m,s} & =\frac{\beta-1}{1-\alpha},\\
\partial_{\theta_{1,n}^{m}}R^{m,s} & =\left(\beta+\frac{1-\beta}{1-\alpha}\right)\sum_{t=1}^{T}\gamma^{t}\left(g_{nt}^{m,s}\left(h_{\theta,nt}^{m,s}-h_{\theta,n\left(t-1\right)}^{m,s}\right)-\mu_{nt}^{m,s}h_{\theta,n\left(t-1\right)}^{m,s}\right), & \forall n\in\mathcal{N},
\end{alignat*}
}and{\small{}
\begin{multline*}
R_{0}^{m,s}=\left(\beta+\frac{1-\beta}{1-\alpha}\right)\sum_{t=1}^{T}\gamma^{t}\Bigg(\sum_{n\in\mathcal{N}}\left(g_{nt}^{m,s}h_{0,nt}^{m,s}-g_{nt}^{m,s}h_{0,n\left(t-1\right)}^{m,s}+g_{0,nt}^{m,s}\right)\\
+\sum_{i\in\mathcal{I}}b_{it}\xi_{it}^{s}-\sum_{i\in\mathcal{I}}\psi_{it}^{m,s}\xi_{it}^{s}\Bigg)+\frac{1-\beta}{1-\alpha}\sum_{n\in\mathcal{N}}g_{0,n0}^{m,s},
\end{multline*}
}{\small\par}
\end{enumerate}
such that
\begin{equation}
R^{s}\left(K_{0}^{m},u^{m},\theta^{m}\right)=\left(\begin{array}{c}
\partial_{K_{0}^{m}}R^{m,s}\\
\partial_{u^{m}}R^{m,s}\\
\partial_{\theta_{1}^{m}}R^{m,s}
\end{array}\right)^{\top}\left(\begin{array}{c}
K_{0}^{m}\\
u^{m}\\
\theta^{m}
\end{array}\right)+R_{0}^{m,s}.\label{eq:gradient_recourse_func}
\end{equation}
\end{prop}
\renewcommand{\labelenumi}{\arabic{enumi}.}

Eq. (\ref{eq:gradient_recourse_func}) indicates that $R^{s}\left(K_{0}^{m},u^{m},\theta^{m}\right)$,
the value of the recourse function at $\left(K_{0}^{m},u^{m},\theta^{m}\right)$,
can be recovered by a linear function with respect to $\left(K_{0}^{m},u^{m},\theta^{m}\right)$
if $\partial R^{m,s}$ and $R_{0}^{m,s}$ are given. Denote 
\[
\phi^{m}\triangleq\mathbb{E}_{s\in\mathcal{S}}\left[\begin{array}{c}
\partial_{K_{0}^{m}}R^{s}\\
\partial_{u^{m}}R^{s}\\
\partial_{\theta_{1}^{m}}R^{s}
\end{array}\right]\ \ \text{and}\ \ \phi_{0}^{m}\triangleq\mathbb{E}_{s\in\mathcal{S}}\left[R_{0}^{m,s}\right],
\]
then a subgradient cut is given by
\begin{equation}
y\ge\left(\phi^{m}\right)^{\top}\left(K_{0},u,\theta_{1}\right)+\phi_{0}^{m}.\label{eq:subgradient_cut}
\end{equation}
We see that if $\left(K_{0},u,\theta_{1}\right)=\left(K_{0}^{m},u^{m},\theta_{1}^{m}\right)$,
the cut recovers $\mathbb{E}\left[R^{m,s}\left(K_{0}^{m},u^{m},\theta^{m}\right)\right]$,
otherwise it returns the recourse along the computed subgradient.
Then, we can update $\left(K_{0},u,\theta_{1}\right)$ by adding the
subgradient cut, Eq. (\ref{eq:subgradient_cut}), into Problem \eqref{master_problem}.

\subsubsection*{Step 3: Update $\left(K_{0},u,\theta_{1}\right)$ by Solving the
First-Stage Problem}

Suppose we have a set of subgradient cuts computed from Step 2 up
to iteration $m$. We add the subgradient cuts to Problem \eqref{master_problem},
and solve the problem in iteration $m+1$. Once solved, we denote
its optimal solution as $\left(K_{0}^{m+1},u^{m+1},\theta_{1}^{m+1}\right)$
and go to Step 1. Problem \eqref{master_problem} is a small-scale
MILP, which has $N+N\times L$ integer variables, and can be efficiently
solved with commercial solvers.

The decomposition algorithm terminates when the objective value of
Problem (\ref{eq:first_stage_problem}), computed from the best-found
solution $\left(K_{0},u,\theta\right)$, is close enough to the objective
value derived from Problem \eqref{master_problem} \emph{or} when
a preset number of iterations is reached. Algorithm \ref{alg:algorithm_procedure}
overviews the details of the entire procedure.

\begin{algorithm}
\caption{\label{alg:algorithm_procedure}Subgradient-based decomposition algorithm}

{\small \begin{algorithmic}[1] \Require{$\bar{m}, \bar{k}, \epsilon$, initial solution $(K_0^1, u^1, \theta_1^1)$}
\Ensure{$(K_0^*, u^*, \theta^*)$}
\State{Initialize $m=1, V_{\text{lb}}= -\infty, V_{\text{ub}}= +\infty$}
\While{$m \le \bar{m}$ $\text{\bf{and}}$ $V_{\text{ub}} - V_{\text{lb}} > \epsilon$ }
\While{$k \le \bar{k}$ \textbf{or} $ \left\lVert \bar{\theta}_2^{m,k+1} - \bar{\theta}_2^{m,k} \right\rVert > \epsilon$}
\State{Randomly pick out a sample path $\{d_1, \ldots, d_T^s\}$ from $\mathcal{S}$}
\State{Solve Problem \eqref{eq:ssp_model} and derive the optimal solution $\hat{\theta}_2^{m,k}$ and $\hat{\delta}_2^{m,k}$}
\State{Update $\bar{\theta}_2^{m,k+1} = \bar{\theta}_2^{m,k} + \sigma_k \left(\hat{\theta}_2^{m,k} - \bar{\theta}_2^{m,k} \right) $}
\State{$k^* \leftarrow k + 1$ and $k \leftarrow k + 1$}
\EndWhile
\State{For all $n \in \mathcal{N}$ and $t \in \mathcal{T}$, set $ \theta_2^m = M_{\theta} \text{ \textbf{if} } \sum_{k^{\prime}=1}^{k^* - 1} \hat{\delta}_{nt}^{m,k^\prime} = 0$; otherwise $\theta_2^m = \bar{\theta}_{2,nt}^{m,k^*} $}
\State{Given $(K_0^m, u^m, \theta^m)$ and Proposition \ref{prop:subgradient}, evaluate $\mathbb{E}[R^s(K_0^m, u^m, \theta^m)]$, and construct Cut \eqref{eq:subgradient_cut}} 			
\State{Solve Problem \eqref{eq:problem_in_step3}-\eqref{eq:end_problem_in_step3} and derive the optimal solution $(K_0^{*}, u^{*}, \theta_1^{*})$ and its optimal value $V^*$}
\State{$V_{\text{ub}} \leftarrow \min \{c_0(K_0^m) + (1-\beta)u^m + \mathbb{E}[R^s(K_0^m, u^m, \theta^m)], V_{\text{ub}} \}$}
\State{$V_{\text{lb}} \leftarrow \max \{V^*, V_{\text{lb}} \}$}
\State{$m \leftarrow m+1$}
\EndWhile 			
\State{$ m^* = \arg \min_{m} \{ c_0(K_0^m) + (1-\beta)u^m + \mathbb{E}[R^s(K_0^m, u^m, \theta^m)] \}$ and $(K_0^*, u^*, \theta^*) = (K_0^{m^*}, u^{m^*}, \theta^{m^*})$ }
\end{algorithmic}}
\end{algorithm}

\subsection{Improving the Best-Found Solution}

The solution $\left(K_{0},u,\theta_{1}\right)$ obtained in Step 3
is not necessarily globally optimal. In Problem \eqref{master_problem},
the epigraph of the recourse function is successively approximated
from below by the subgradient cuts, i.e. Eq. (\ref{eq:subgradient_cut}).
If the recourse function is convex, this procedure is the same as
Benders decomposition, and the global optimum can be found since the
subgradient gives a global underestimator for convex recourse functions
\citep{benders_partitioning_1962}. However, this is not the case
for $R^{s}\left(K_{0},u,\theta\right)$ which is non-convex due to
the integer wait-and-see variables and the non-convex cost function
$c_{t}\left(\cdot\right)$. In this case, the subgradient cuts may
cut off part of the epigraph. Therefore, we wish to further improve
upon the best-found solution.

We use the multi-cut method proposed in \citep{zhao_approximate_2019}
to improve the solution given by Algorithm \ref{alg:algorithm_procedure}.
First, we introduce another type of cut, i.e. the integer optimality
cut \citep{laporte_integer_1993}, simultaneously with the subgradient
cuts in Problem \eqref{master_problem} to update $\left(K_{0},u,\theta\right)$.
The integer optimality cut is valid, i.e. it is a global underestimator
of the epigraph, but it is conservative when updating the first-stage
decisions. Second, we run the decomposition algorithm until the termination
condition is triggered. Then, we stop adding subgradient cuts and
remove one of the subgradient cuts which has the minimum slack. This
procedure relaxes the approximation of the epigraph and helps the
algorithm avoid getting stuck at a local solution.

\section{\label{sec:Numerical-Study}Case Study: Capacity Planning for a Waste-to-Energy
System}

Our case study about a multi-facility WTE system in Singapore is adapted
from \citep{cardin_analyzing_2015,zhao_decision_2018}. All data in
this section are taken from the real case study.

The system has five candidate sites in different parts of Singapore.
The WTE facility disposes of food waste collected from each sector
by using an innovative anaerobic digestion technique which transforms
the food waste into electricity. Undisposed waste will be subjected
to further treatment via landfill, incurring greater disposal costs
(i.e. penalties). The profit of the system comes from selling the
electricity, and the costs consist of unit disposal costs, penalty
costs, transportation costs, and capacity expansion costs. For simplicity,
we omit annual fixed costs incurred once the facility is installed,
which differs from the models in \citep{cardin_analyzing_2015,zhao_decision_2018}.
However, the proposed method can still solve the problem with fixed
costs if we use a multi-cut version, i.e. subgradient cuts with integer
optimality cuts, in Algorithm \ref{alg:algorithm_procedure}.

The generation of the food waste from each sector (i.e. the demand
of each customer) is assumed to be standard geometric Brownian motion
(GBM) \citep{cardin_analyzing_2015}, represented by
\[
\xi_{it}=\left(\bar{\mu}+\bar{\sigma}\omega_{t}\right)\xi_{i\left(t-1\right)},\ \ \forall t\in\mathcal{T},
\]
where $\xi_{it}$ is the waste amount generated from sector $i$ in
time $t$, $\bar{\mu}$ is the percentage drift, $\bar{\sigma}$ is
the percentage volatility, and $\omega_{t}$ is a standard normal
random variable. In the numerical study, we assume that $\bar{\mu}$
is 4\% and $\bar{\sigma}$ is 16\%, and the initial waste vector at
the beginning is $\xi_{0}=\left[498,518,293,460,382\right]$ (unit:
tonnes per day).

The WTE facilities are assumed to be modular; that is, one unit of
capacity can dispose of 100 tonnes of food waste per day, and the
capacity limit is $K^{\text{max}}=\left[16,10,10,10,10\right]$. The
expansion cost of facility $n$ is given by a power function:
\[
\bar{c}_{nt}\left(\Delta K_{nt}\right)=7\times10^{6}\times\left(\Delta K_{nt}\right)^{0.9}\times100,\ \ \forall n\in\mathcal{N},t\in\mathcal{T}.
\]
We linearize the expansion costs by setting breakpoints at $\left(a_{1},a_{2},\ldots,a_{7}\right)=\left(0,2^{0},2^{1},\ldots,2^{4}\right)$
and derive the cost function presented in Eq. (\ref{eq:linearize_cost_func}).
The penalty for undisposed food waste is $b_{it}=77$ per tonne, and
the profit matrix is{\small{}
\[
\left(r_{in}\right)_{i\in\mathcal{I},n\in\mathcal{N}}=\left[\begin{array}{ccccc}
59.9 & 42.9 & 30.4 & 34.9 & 21.9\\
42.9 & 59.9 & 40.9 & 52.1 & 38.9\\
30.4 & 40.9 & 59.9 & 49.8 & 34.4\\
34.9 & 52.1 & 49.9 & 59.9 & 44.4\\
21.9 & 38.9 & 34.4 & 44.4 & 59.9
\end{array}\right],
\]
}where $r_{in}$ denotes the net profit from disposing every tonne
of food waste from sector $i$ by facility $n$. Since the system
has a one-to-one correspondence between a sector and its most profitable
facility, the allocation matrix $W$ is a $5\times5$ identity matrix.
Therefore, the decision rule is
\[
K_{nt}\left(\theta\right)=\begin{cases}
\left\lfloor \xi_{it}/100\right\rceil +\theta_{1,n}, & \text{if }\left\lfloor \xi_{it}/100\right\rceil -K_{n\left(t-1\right)}\ge\theta_{2,nt},\\
K_{n\left(t-1\right)}, & \text{otherwise.}
\end{cases}
\]
We set the discount factor $\gamma=0.926$.

\subsection{A Baseline Design and the Value of Flexibility}

Problem (\ref{eq:Risk_aversion_Model}) is a multi-stage stochastic
programming problem and it is very difficult to solve. The solution
provided by the decision-rule based method may not be globally optimal.
Therefore, we need an attainable baseline design as a lower bound,
to see how the proposed method improves upon the economic performance.

A baseline is a restriction of Problem (\ref{eq:Risk_Neutral_Model}):
the decision maker fixes the capacity plan at the beginning, and has
no options to adjust the capacity in the future. Mathematically, the
capacity decisions $\left(K_{0},K_{1},\ldots,K_{T}\right)$ of the
baseline model, once determined, do not change with the realizations
of $\xi$ over time. We call this baseline MCEP the \emph{inflexible}
model. This inflexible model can be solved by the Benders decomposition:
we fix the capacity decisions, i.e. $\left(K_{0},K_{1},\ldots,K_{T}\right)$,
solve the allocation problems (i.e. $\Pi_{t}\left(\cdot\right)$ for
all $t\in\mathcal{T}$), and then iteratively update $\left(K_{0},K_{1},\ldots,K_{T}\right)$
via the generated Benders cuts. The allocation problems are all LPs,
and we have finite number of capacity decisions, so the optimal solution
of the inflexible model will be reached within a finite number of
iterations \citep{benders_partitioning_1962}.

The flexible and inflexible models are solved separately by generating
4,000 demand scenarios via Monte-Carlo simulation. Then, to compare,
we conduct out-of-sample tests to evaluate the economic performance
of each design by using an identical sample set. The out-of-sample
test consists of 12,000 scenarios that are generated via Monte-Carlo
simulation. The rationale underlying the out-of-sample test is to
eliminate the bias introduced by using different sample sets.

Given the baseline design, we can calculate the difference between
the economic performance of the flexible multi-stage model and its
inflexible counterpart, i.e. the value of flexibility (VoF).

\subsection{The Proposed Method Captures Decision-Maker Risk-Preferences}

In our simulations, we set $\alpha=0.95$ and vary $\beta$ from $0.01$
to $0.99$. When $\alpha=0.95$, we are minimizing expected costs
that are greater than or equal to the $95^{\text{th}}$-percentile.
Equivalently, the objective of Problem (\ref{eq:Risk_aversion_Model})
maximizes the expected value less than the $5^{\text{th}}$-percentile.
The simulation results from different policies are presented in Table
\ref{tab:SimulationResults}. All cases are run three times and the
values are averaged. The inflexible policy comes from the baseline
design, where ENPV is 237.7 million, the $5^{\text{th}}$-percentile
of the NPV is 146.4 million, and the worst-case NPV is $-291.8$ million.
The flexible policies from Problem (\ref{eq:Risk_aversion_Model})
are computed for varying $\beta$. If we set $\beta=0.99$, the ENPV
is 291.3 million, the $5^{\text{th}}$-percentile is 217.1 million,
and the worst-case scenario is 130.1 million. We see that the flexible
policy dominates the inflexible one in terms of the five metrics presented
in Table \ref{tab:SimulationResults}. In particular, the worst-case
scenario is improved significantly. The worst-case NPV for the inflexible
policy is negative, while the worst-case NPV of the flexible policies
ranges from 125.1 million to 130.1 million.

\begin{table}
\begin{doublespace}
\caption{\label{tab:SimulationResults}Simulation results with $\alpha=0.95$
(unit: million S\$).}

\end{doublespace}
\centering{}%
\begin{tabular}{ccccccc}
\toprule 
\multirow{2}{*}{Criterion} & \multicolumn{5}{c}{Flexible policies with different $\beta$} & \multirow{2}{*}{Inflexible policy}\tabularnewline
\cmidrule{2-6} \cmidrule{3-6} \cmidrule{4-6} \cmidrule{5-6} \cmidrule{6-6} 
 & 0.01 & 0.25 & 0.5 & 0.75 & 0.99 & \tabularnewline
\midrule
Min & 125.1 & 135.5 & 135.0 & 135.5 & 130.1 & $-291.8$\tabularnewline
$5^{\text{th}}$-percentile & 224.3 & 221.4 & 220.5 & 218.7 & 217.1 & 146.4\tabularnewline
ENPV & 284.8 & 287.4 & 288.7 & 290.3 & 291.3 & 237.7\tabularnewline
$95^{\text{th}}$-percentile & 344.3 & 354.3 & 361.1 & 365.1 & 366.2 & 306.4\tabularnewline
Max & 492.5 & 516.6 & 523.0 & 539.9 & 542.3 & 342.4\tabularnewline
\midrule 
VoF & 47.1 & 49.7 & 51.0 & 52.6 & 53.6 & -\tabularnewline
\bottomrule
\end{tabular}
\end{table}

\begin{figure}
\begin{centering}
\includegraphics[scale=0.63]{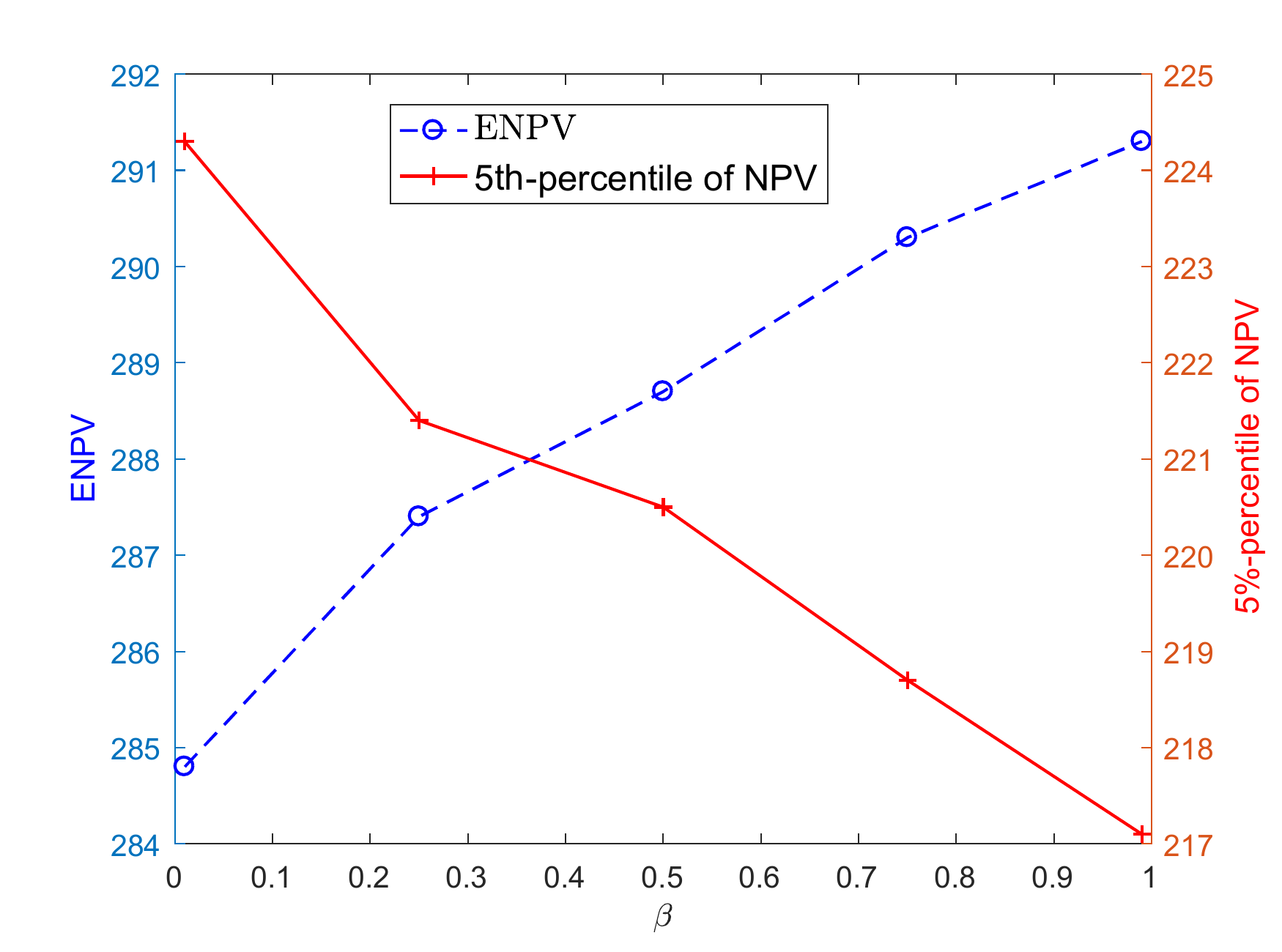}
\par\end{centering}
\caption{\label{fig:Sen_A_on_beta}Simulation results given different $\beta$.}
\end{figure}

These results demonstrate that the $5^{\text{th}}$-percentile of
the NPVs decreases with $\beta$ and the ENPV increases with $\beta$
(see Figure \ref{fig:Sen_A_on_beta}). If we change $\beta$ from
0.99 to 0.01, the $5^{\text{th}}$-percentile rises from 217.1 million
to 224.3 million and the ENPV declines from 291.3 million to 284.8
million. In addition, the VoF decreases from 53.6 million to 47.1
million, which means that the value gained from flexibility decreases
as the weight factor decreases. One possible explanation is because
the particular flexible capacity-expansion policy analyzed here focuses
on improving upside potential. This means that risk-averse decision-makers
may gain less value from being flexible. The result from Figure \ref{fig:Sen_A_on_beta}
also verifies the conclusion of Proposition \ref{prop:ENPV_propertise}
numerically.

\subsection{Risk-Averse Policy is More Conservative in Expansion}

\begin{figure}
\begin{centering}
\includegraphics[viewport=25bp 0bp 680bp 360bp,clip,scale=0.63]{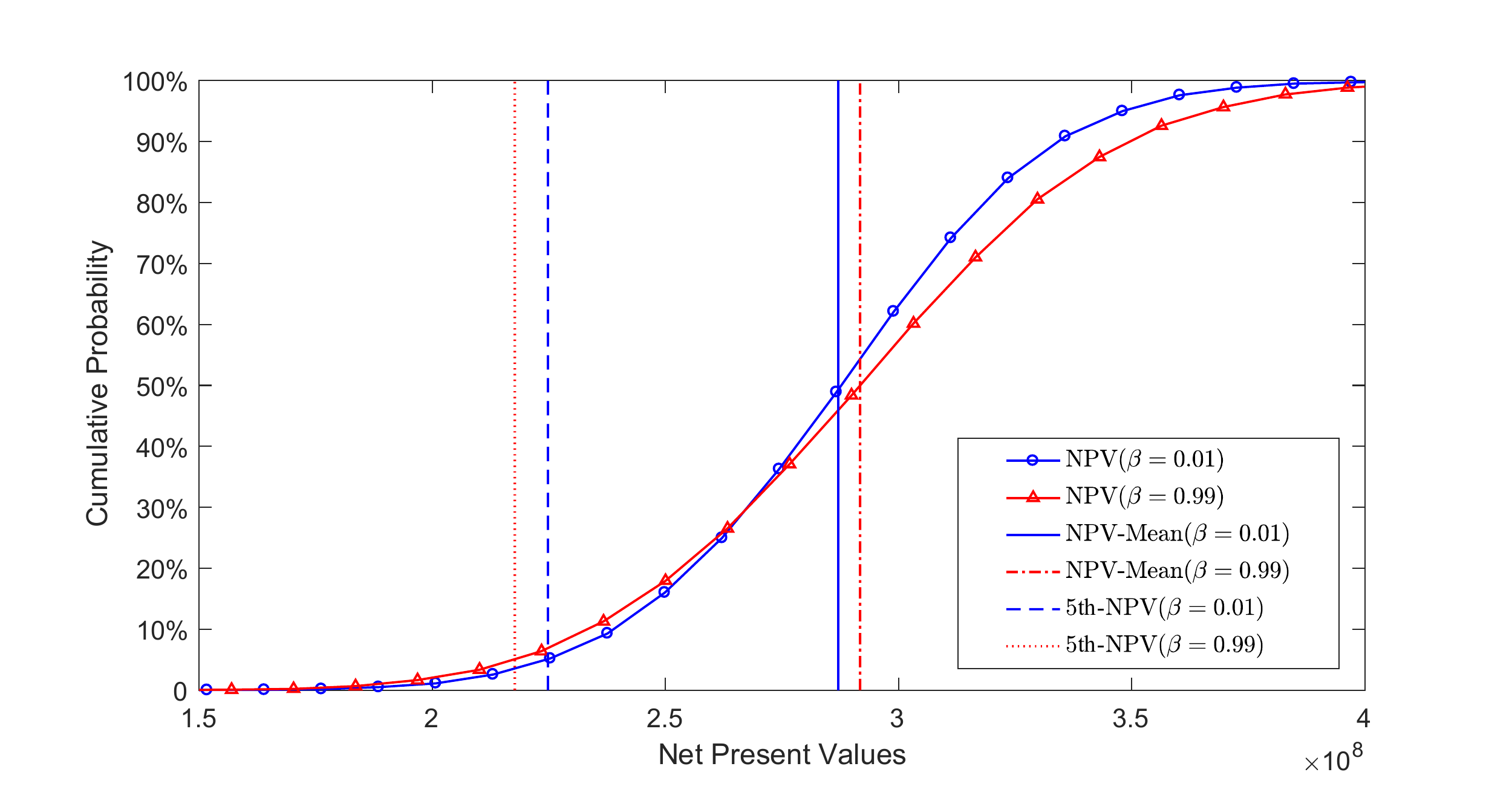}
\par\end{centering}
\caption{\label{fig:CDF_comparison}Cumulative distribution functions of different
policies.}
\end{figure}

We examine the optimal solutions of Problem (\ref{eq:Risk_aversion_Model})
for different risk preferences. The optimal initial capacity given
$\beta=0.99$ is $K_{0}^{*}=\left[6,7,3,6,4\right]$ and the optimal
parameter is $\theta_{1}^{*}=\left[0,0,0,1,0\right]$. In contrast,
the optimal solution given $\beta=0.01$ is $K_{0}^{*}=\left[5,6,3,5,4\right]$
and $\theta_{1}^{*}=\left[0,0,0,0,0\right]$. The policy with $\beta=0.01$
has a smaller initial capacity, and its expansion is more conservative
when the decision rule is triggered. The cumulative density functions
of the out-of-sample tests of these two policies are plotted in Figure
\ref{fig:CDF_comparison}. As can be seen, the policy with $\beta=0.01$
does not fully exploit the upside expansion opportunity whereas the
one with $\beta=0.99$ does, but it reduces the downside cost as it
is more conservative.

\subsection{Flexible Policies have Better Performance under an Inaccurate Demand
Model}

We also test the robustness of the policies when the demand model
is inaccurate. Recall that our demand model is GBM with $\bar{\mu}=0.04$
and $\bar{\sigma}=0.16$. We perform out-of-sample tests via a set
of samples generated from a GBM with \emph{identical} $\bar{\mu}$
and $\bar{\sigma}$. To see how the performance changes if we use
incorrect training data, we intentionally generate samples via GBMs
that have \emph{different} parameters and perform out-of-sample tests.
We generate two sample sets:
\begin{itemize}
\item Set A with $12,000$ samples generated via $\text{GBM}\left(0.02,0.08\right)$,
which has a lower percentage drift and lower percentage volatility
compared to the training data;
\item Set B with $12,000$ samples generated via $\text{GBM}\left(0.05,0.16\right)$,
which has a higher percentage drift and higher percentage volatility
compared to the training data.
\end{itemize}
The simulation results of the out-of-sample tests are presented in
Table \ref{tab:incorrect_prediction}. We see that for the risk-averse
policy (with $\beta=0.01$), the ENPV computed under the original
sampling condition is $284.8$ million, but if we perform out-of-sample
tests via the alternative sample Sets A and B, its ENPV decreases
$2.3\%$ and increases $3.1\%$, respectively. The ENPV derived from
the risk-neutral flexible policy (with $\beta=0.99$) decreases $4.2\%$
in Set A and increases $5.5\%$ in Set B. However, the ENPVs of the
inflexible policy decrease $-10.3\%$ in Set A and increase merely
$1.8\%$ in Set B. Several conclusions follow from Table \ref{tab:incorrect_prediction}.
First, a risk-averse policy may suffer lower loss compared to a risk-neutral
policy when the predicted demands (i.e., the training data) are not
as large as the true demands (i.e., the out-of-test data). Second,
we see that the risk-neutral policy can better exploit upside expansion
opportunities if the testing demands are generally larger than the
training demands. Third, flexible policies are more robust in comparison
with inflexible policies for an incorrect demand model.

\begin{table}
\caption{\label{tab:incorrect_prediction}Simulation results given inaccurate
demand models (unit: million S\$).}

\centering{}%
\begin{tabular}{ccccccc}
\toprule 
\multirow{2}{*}{Policies} & \multirow{2}{*}{Criteria} & \multirow{2}{*}{$\begin{array}{c}
\text{Training samples}\\
\text{GBM\ensuremath{\left(0.04,0.16\right)}}
\end{array}$} & \multicolumn{2}{c}{Set A} & \multicolumn{2}{c}{Set B}\tabularnewline
\cmidrule{4-7} \cmidrule{5-7} \cmidrule{6-7} \cmidrule{7-7} 
 &  &  & Results & Variations & Results & Variations\tabularnewline
\midrule
\multirow{2}{*}{$\begin{array}{c}
\text{Flexible policy}\\
(\beta=0.01)
\end{array}$} & $5^{\text{th}}$-percentile & 224.3 & 249.4 & $+11.2\%$ & 236.4 & $+5.1\%$\tabularnewline
 & ENPV & 284.8 & 278.2 & $-2.3\%$ & 295.9 & $+3.1\%$\tabularnewline
\midrule
\multirow{2}{*}{$\begin{array}{c}
\text{Flexible policy}\\
(\beta=0.99)
\end{array}$} & $5^{\text{th}}$-percentile & 217.1 & 239.7 & $+10.4\%$ & 233.2 & $+7.1\%$\tabularnewline
 & ENPV & 291.3 & 279.2 & $-4.2\%$ & 307.9 & $+5.5\%$\tabularnewline
\midrule
\multirow{2}{*}{Inflexible policy} & $5^{\text{th}}$-percentile & 146.4 & 157.5 & $+7.6\%$ & 147.5 & $+0.8\%$\tabularnewline
 & ENPV & 237.7 & 213.2 & $-10.3\%$ & 242.2 & $+1.8\%$\tabularnewline
\bottomrule
\end{tabular}
\end{table}

\section{Conclusion}

In this paper, we establish a multi-facility capacity expansion model
with a mean-CVaR objective to capture decision maker risk preferences.
We verify that the economic performance of the risk-averse model is
not higher than that of the risk-neutral model.

To solve the our risk-averse model, we approximate the expansion policy
of the multi-stage problem with an if-then decision rule, and then
solve the resulting model by a customized subgradient decomposition
algorithm. Though our algorithm may not reach the global optimum of
the problem, numerical studies show that its improvement over the
baseline design (i.e. the inflexible MCEP) is significant, especially
when the demand model is inaccurate.

Simulation results also reveal that the decision maker is able to
choose a policy with a higher ENPV or one with a higher $5^{\text{th}}$-percentile
of the NPVs by simply adjusting the weight factor of the objective
function. In addition, the ENPV decreases as the decision maker becomes
more risk-averse, indicating that the decision maker will prefer to
pay less for flexibility. The simulation results also illustrate that
a risk-averse expansion policy, compared with the risk-neutral one,
tends to establish smaller initial capacity at the beginning, and
is more conservative in future expansion.

Many opportunities exist for future work. The model proposed in this
paper does not account for annual fixed costs. Such costs can take
the form of annual land rental costs which are incurred once a facility
is established (i.e. has a non-zero capacity). To formulate a problem
with fixed costs, we need to introduce extra binary auxiliary variables
as the costs are concave with respect to the capacity, so the model
becomes more difficult to handle numerically. In another research
direction, we may consider robustness against an unknown demand model.
Our present numerical experiments suggest that our approach has some
intrinsic robustness against demand uncertainty, so we can build on
this initial proof of concept.

\appendix

\section*{Appendix. \label{sec:Appendix-A}Proofs of Main Results}

\noun{Proof of Proposition \ref{prop:ENPV_propertise}:} Firstly,
we know from the definition of CVaR that, given a random variable
$X$, we have $\text{CVaR}_{\alpha}\left(X\right)\ge\mathbb{E}\left[X\right]$
for any $\alpha\in\left(0,1\right)$ \citep{chen_value-at-risk_2008}.

(i) Denote $\mathcal{K}_{\left[T\right]}^{*}$ as the optimal solution
of Problem (\ref{eq:Risk_aversion_Model}) given $\beta$. Denote
$\Delta\beta>0$ as a positive increment such that $\beta+\Delta\beta\le1$.
According to Problem (\ref{eq:Risk_aversion_Model}), we have
\begin{alignat*}{1}
\text{ENPV}_{\alpha}\left(\beta+\Delta\beta\right) & =\max_{K_{\left[T\right]}\in\mathcal{K}}\Big(-\beta\mathbb{E}\left[Q\left(\mathcal{K}_{\left[T\right]},\xi\right)\right]-\left(1-\beta\right)\text{CVaR}_{\alpha}\left(Q\left(\mathcal{K}_{\left[T\right]},\xi\right)\right)\\
 & \ \ \ \ \ \ \ \ \ \ \ \ \ \ \ \ \ \ \ \ \ \ \ \ \ \ \ \ \ \ \ \ \ \ +\Delta\beta\left(\text{CVaR}_{\alpha}\left(Q\left(\mathcal{K}_{\left[T\right]},\xi\right)\right)-\mathbb{E}\left[Q\left(\mathcal{K}_{\left[T\right]},\xi\right)\right]\right)\Big)\\
 & \ge-\beta\mathbb{E}\left[Q\left(\mathcal{K}_{\left[T\right]}^{*},\xi\right)\right]-\left(1-\beta\right)\text{CVaR}_{\alpha}\left(Q\left(\mathcal{K}_{\left[T\right]}^{*},\xi\right)\right)\\
 & \ \ \ \ \ \ \ \ \ \ \ \ \ \ \ \ \ \ \ \ \ \ \ \ \ \ \ \ \ \ \ \ \ \ +\Delta\beta\left(\text{CVaR}_{\alpha}\left(Q\left(\mathcal{K}_{\left[T\right]}^{*},\xi\right)\right)-\mathbb{E}\left[Q\left(\mathcal{K}_{\left[T\right]}^{*},\xi\right)\right]\right)\\
 & =\text{ENPV}_{\alpha}\left(\beta\right)+\Delta\beta\left(\text{CVaR}_{\alpha}\left(Q\left(\mathcal{K}_{\left[T\right]}^{*},\xi\right)\right)-\mathbb{E}\left[Q\left(\mathcal{K}_{\left[T\right]}^{*},\xi\right)\right]\right)
\end{alignat*}
The inequality holds as $\mathcal{K}_{\left[T\right]}^{*}$ is the
optimal solution for $\text{ENPV}_{\alpha}\left(\beta\right)$. Since
$\text{CVaR}_{\alpha}\left(Q\left(\mathcal{K}_{\left[T\right]},\xi\right)\right)-\mathbb{E}\left[Q\left(\mathcal{K}_{\left[T\right]},\xi\right)\right]\ge0$
for any $\mathcal{K}_{\left[T\right]}\in\bar{\mathcal{K}}$, we can
conclude that $\text{ENPV}_{\alpha}\left(\beta\right)$ is non-decreasing
in $\beta$.

(ii) Since $\text{ENPV}_{\alpha}\left(1\right)=\text{ENPV}_{\text{flex}}$
and $\text{ENPV}_{\alpha}\left(\beta\right)$ is non-decreasing in
$\beta$, given any $\alpha\in\left(0,1\right)$, we have $\text{ENPV}_{\text{inflex}}=\text{ENPV}_{\alpha}\left(1\right)\ge\text{ENPV}_{\alpha}\left(\beta\right)$
for all $\beta\in\left[0,1\right]$ according to (i). Hence, we can
conclude the result.\qed

 \ \ 

\noindent\noun{Proof of Lemma \ref{lem:closed_form_for_K}:} We first
show that {\small{}
\[
K_{nt}^{m,s}=\left[\prod_{t_{0}=1}^{t}\left(1-\delta_{nt_{0}}^{m,s}\right)\right]K_{n0}^{m}+\sum_{t_{0}=1}^{t}\left[\prod_{t_{1}=t_{0}+1}^{t}\left(1-\delta_{nt_{1}}^{m,s}\right)\right]\delta_{nt_{0}}^{m,s}\left(\left\lfloor \sum\limits _{i\in\mathcal{I}}W_{in}\xi_{it_{0}}^{s}\right\rceil +\theta_{1,n}^{m}\right),\ \ \forall n\in\mathcal{N},t\in\mathcal{T}\cup\left\{ 0\right\} .
\]
}The result can be proved by induction. For $t=1$, we have 
\begin{alignat*}{1}
K_{n1}^{m,s} & =\left(1-\delta_{n1}^{m,s}\right)K_{n0}^{m}+\delta_{n1}^{m,s}\left(\left\lfloor \sum\limits _{i\in\mathcal{I}}W_{in}\xi_{i1}^{s}\right\rceil +\theta_{1,n}^{m}\right),
\end{alignat*}
and the result holds. Suppose the result holds for $t=t^{\prime}-1$:
\[
K_{n\left(t^{\prime}-1\right)}^{m,s}=\left[\prod_{t_{0}=1}^{t^{\prime}-1}\left(1-\delta_{nt_{0}}^{m,s}\right)\right]K_{n0}^{m}+\sum_{t_{0}=1}^{t^{\prime}-1}\left[\prod_{t_{1}=t_{0}+1}^{t^{\prime}-1}\left(1-\delta_{nt_{1}}^{m,s}\right)\right]\delta_{nt_{0}}^{m,s}\left(\left\lfloor \sum\limits _{i\in\mathcal{I}}W_{in}\xi_{it_{0}}^{s}\right\rceil +\theta_{1,n}^{m}\right).
\]
Then, for $t=t^{\prime}$ we have
\begin{alignat*}{1}
K_{nt^{\prime}}^{m,s} & =\left(1-\delta_{nt^{\prime}}^{m,s}\right)K_{n\left(t^{\prime}-1\right)}^{m,s}+\delta_{nt^{\prime}}^{m,s}\left(\left\lfloor \sum\limits _{i\in\mathcal{I}}W_{in}\xi_{it^{\prime}}^{s}\right\rceil +\theta_{1,n}^{m}\right)\\
 & =\left[\prod_{t_{0}=1}^{t^{\prime}}\left(1-\delta_{nt_{0}}^{m,s}\right)\right]K_{n0}^{m}\\
 & \ \ \ \ \ \ \ \ +\sum_{t_{0}=1}^{t^{\prime}-1}\left[\prod_{t_{1}=t_{0}+1}^{t^{\prime}}\left(1-\delta_{nt_{1}}^{m,s}\right)\right]\delta_{nt_{0}}^{m,s}\left(\left\lfloor \sum\limits _{i\in\mathcal{I}}W_{in}\xi_{it_{0}}^{s}\right\rceil +\theta_{1,n}^{m}\right)+\delta_{nt^{\prime}}^{m,s}\left(\left\lfloor \sum\limits _{i\in\mathcal{I}}W_{in}\xi_{it^{\prime}}^{s}\right\rceil +\theta_{1,n}^{m}\right)\\
 & =\left[\prod_{t_{0}=1}^{t^{\prime}}\left(1-\delta_{nt_{0}}^{m,s}\right)\right]K_{n0}^{m}+\sum_{t_{0}=1}^{t^{\prime}}\left[\prod_{t_{1}=t_{0}+1}^{t^{\prime}}\left(1-\delta_{nt_{1}}^{m,s}\right)\right]\delta_{nt_{0}}^{m,s}\left(\left\lfloor \sum\limits _{i\in\mathcal{I}}W_{in}\xi_{it_{0}}^{s}\right\rceil +\theta_{1,n}^{m}\right).
\end{alignat*}
The last equality holds as we have $\prod_{t_{1}=t_{0}+1}^{t^{\prime}}\left(1-\delta_{nt_{1}}^{m,s}\right)=\prod_{t_{1}=t^{\prime}+1}^{t^{\prime}}\left(1-\delta_{nt_{1}}^{m,s}\right)=1$.
Reorganize the above equation and we can conclude the result.\qed

 \ \ 

\noindent\noun{Proof of Lemma \ref{lem:grad_C}: }According to Eq.
(\ref{eq:linearize_cost_func}) and the definitions of $g_{nt}^{m,s}$,
we have{\small{}
\begin{alignat*}{1}
\partial c_{t}^{m,s}= & \ \partial\left(\sum_{n\in\mathcal{N}}\left[g_{nt}^{m,s}\cdot\left(K_{nt}^{m,s}-K_{n\left(t-1\right)}^{m,s}\right)+g_{0,nt}^{m,s}\right]\right),\\
= & \ \partial\left(\sum_{n\in\mathcal{N}}\left[g_{nt}^{m,s}\cdot\left(h_{K,nt}^{m,s}K_{n0}^{m}+h_{\theta,nt}^{m,s}\theta_{1,n}^{m}+h_{0,nt}^{m,s}-h_{K,n\left(t-1\right)}^{m,s}K_{n0}^{m}-h_{\theta,n\left(t-1\right)}^{m,s}\theta_{1,n}^{m}-h_{0,n\left(t-1\right)}^{m,s}\right)+g_{0,nt}^{m,s}\right]\right).
\end{alignat*}
}Then, the following result holds trivially according to Eq. (\ref{eq:closed_form_for_K})
from Lemma \ref{lem:closed_form_for_K}. \qed

 \ \ 

\noindent\noun{Proof of Lemma \ref{lem:grad_Pi}:} Let $\left(\mu_{nt}\right)_{n\in\mathcal{N},t\in\mathcal{T}}$
and $\left(\psi_{it}\right)_{i\in\mathcal{I},t\in\mathcal{T}}$ be
the dual variables with respect to the capacity and demand constraints
of $\Pi_{t}\left(\cdot\right)$ respectively, and the dual problem
of $\Pi_{t}\left(K_{t-1}^{m,s},\xi_{t}^{s}\right)$ can then be formulated
by
\begin{alignat}{2}
\min_{\psi,\mu}\  & \sum_{n\in\mathcal{N}}\mu_{nt}K_{n\left(t-1\right)}^{m,s}+\sum_{i\in\mathcal{I}}\psi_{it}\xi_{it}^{s}-\sum_{i\in\mathcal{I}}b_{it}\xi_{it}^{s}\label{eq:dual_of_pi}\\
\text{s.t.\  } & -r_{int}-b_{it}+\psi_{it}+\mu_{nt}\ge0, &  & \ \ \forall i\in\mathcal{I},n\in\mathcal{N},\nonumber \\
 & \psi_{it}\ge0,\ \ \mu_{nt}\ge0, &  & \ \ \forall i\in\mathcal{I},n\in\mathcal{N}.\nonumber 
\end{alignat}
\textit{\emph{Since Problem (\ref{eq:Allocation_model}) is a linear
program and the optimal solution is not empty as $K_{t-1}^{m,s}$
and $\xi_{t}^{s}$ are all bounded and nonnegative, the strong duality
holds for Problem (\ref{eq:Allocation_model}) \citep[Theorem 4.4]{bertsimas1997introduction}.
Therefore, we have
\[
\Pi_{t}\left(K_{t-1}^{m,s},\xi_{t}^{s}\right)=\sum_{n\in\mathcal{N}}\mu_{nt}^{m,s}K_{n\left(t-1\right)}^{m,s}+\sum_{i\in\mathcal{I}}\psi_{it}^{m,s}\xi_{it}^{s}-\sum_{i\in\mathcal{I}}b_{it}\xi_{it}^{s},
\]
where $\left(\mu_{nt}^{m,s}\right)_{n\in\mathcal{N},t\in\mathcal{T}}$
and $\left(\psi_{it}^{m,s}\right)_{i\in\mathcal{I},t\in\mathcal{T}}$
are the optimal solutions of Problem (\ref{eq:dual_of_pi}) by definitions.
Then, according to Lemma \ref{lem:closed_form_for_K}, Lemma }}\ref{lem:grad_Pi}
\textit{\emph{can be proved.}} \qed

 \ \ 

\noindent\noun{Proof of Proposition \ref{prop:subgradient}: }We
know that $K_{1}^{s},\ldots,K_{T}^{s},\eta^{s}$ can be computed from
Policy (\ref{eq:Decision_Rule}) once $\left(K_{0}^{m},u^{m},\theta^{m}\right)$
and scenario $s\in\mathcal{S}$ are given. According to Eq. (\ref{eq:second_stage_problem}),
we have
\begin{alignat*}{1}
\partial R^{m,s} & =\partial\left[\beta\sum_{t=1}^{T}\gamma^{t}\left(c_{t}\left(K_{t}^{s}-K_{t-1}^{s}\right)-\Pi_{t}\left(K_{t-1}^{s},\xi_{t}^{s}\right)\right)+\frac{1-\beta}{1-\alpha}\eta^{s}\right],\\
 & =\beta\sum_{t=1}^{T}\gamma^{t}\left(\partial c_{t}^{m,s}-\partial\Pi_{t}^{m,s}\right)+\frac{1-\beta}{1-\alpha}\partial\eta^{m,s},
\end{alignat*}
where $\partial c_{t}^{m,s}$, $\partial\Pi_{t}^{m,s}$, and $\partial\eta^{m,s}$
denote the subgradients of the corresponding functions with respect
to $\left(K_{0}^{m},u^{m},\theta_{1}^{m}\right)$. Then, when $c_{0}\left(K_{0}\right)+\sum_{t=1}^{T}\gamma^{t}\left(c_{t}\left(K_{t}^{s}-K_{t-1}^{s}\right)-\Pi_{t}\left(K_{t-1}^{s},\xi_{t}^{s}\right)\right)-u\ge0$,
we have $\partial\eta^{m,s}=\partial c_{0}^{m,s}+\sum_{t=1}^{T}\gamma^{t}\left(\partial c_{t}^{m,s}-\partial\Pi_{t}^{m,s}\right)+\partial u$;
therefore
\begin{equation}
\partial R^{m,s}=\left(\beta+\frac{1-\beta}{1-\alpha}\right)\sum_{t=1}^{T}\gamma^{t}\left(\partial c_{t}^{m,s}-\partial\Pi_{t}^{m,s}\right)+\frac{1-\beta}{1-\alpha}\left(\partial c_{0}^{m,s}+\partial u\right).\label{eq:A.R1}
\end{equation}
When $c_{0}\left(K_{0}\right)+\sum_{t=1}^{T}\gamma^{t}\left(c_{t}\left(K_{t}^{s}-K_{t-1}^{s}\right)-\Pi_{t}\left(K_{t-1}^{s},\xi_{t}^{s}\right)\right)-u<0$,
we have $\partial\eta^{m,s}=0$, and thus 
\begin{alignat}{1}
\partial R^{m,s} & =\beta\sum_{t=1}^{T}\gamma^{t}\left(\partial c_{t}^{m,s}-\partial\Pi_{t}^{m,s}\right).\label{eq:A.R2}
\end{alignat}
Proposition \ref{prop:subgradient} can then be concluded by substituting
the results from Lemma \ref{lem:grad_C} and \ref{lem:grad_Pi} into
Eqs. (\ref{eq:A.R1})--(\ref{eq:A.R2}). \qed

\bibliographystyle{apa}
\bibliography{MyLibrary}

\end{document}